\documentclass[11pt]{article}

\usepackage[margin=1.05in]{geometry}
\usepackage[T1]{fontenc}
\usepackage[utf8]{inputenc}
\usepackage{lmodern}
\usepackage{amsmath,amssymb,amsthm,mathtools}
\usepackage{aliascnt}
\usepackage{mathrsfs}
\usepackage{microtype}
\usepackage[colorlinks=true,linkcolor=blue,citecolor=blue,urlcolor=blue]{hyperref}
\usepackage[nameinlink,noabbrev]{cleveref}

\numberwithin{equation}{section}

\theoremstyle{plain}
\newtheorem{theorem}{Theorem}[section]

\newaliascnt{proposition}{theorem}
\newtheorem{proposition}[proposition]{Proposition}
\aliascntresetthe{proposition}

\newaliascnt{lemma}{theorem}
\newtheorem{lemma}[lemma]{Lemma}
\aliascntresetthe{lemma}

\newaliascnt{corollary}{theorem}
\newtheorem{corollary}[corollary]{Corollary}
\aliascntresetthe{corollary}

\theoremstyle{definition}
\newaliascnt{assumption}{theorem}
\newtheorem{assumption}[assumption]{Assumption}
\aliascntresetthe{assumption}

\theoremstyle{remark}
\newaliascnt{remark}{theorem}
\newtheorem{remark}[remark]{Remark}
\aliascntresetthe{remark}

\Crefname{theorem}{Theorem}{Theorems}
\Crefname{proposition}{Proposition}{Propositions}
\Crefname{lemma}{Lemma}{Lemmas}
\Crefname{corollary}{Corollary}{Corollaries}
\Crefname{assumption}{Assumption}{Assumptions}
\Crefname{remark}{Remark}{Remarks}

\DeclareMathOperator{\Tr}{Tr}
\DeclareMathOperator{\rank}{rank}
\DeclareMathOperator{\Ad}{Ad}
\DeclareMathOperator{\ad}{ad}
\DeclareMathOperator{\Herm}{Herm}
\DeclareMathOperator{\diag}{diag}
\DeclareMathOperator{\pr}{pr}
\DeclareMathOperator{\Lie}{Lie}

\newcommand{\R}{\mathbb R}
\newcommand{\C}{\mathbb C}
\newcommand{\Z}{\mathbb Z}
\newcommand{\gLie}{\mathfrak g}
\newcommand{\kLie}{\mathfrak k}
\newcommand{\lLie}{\mathfrak l}
\newcommand{\pLie}{\mathfrak p}
\newcommand{\sr}{\mathrm{SR}}
\newcommand{\SU}{\mathrm{SU}}
\newcommand{\U}{\mathrm{U}}
\newcommand{\norm}[1]{\left\lVert #1\right\rVert}
\newcommand{\ip}[2]{\left\langle #1,#2\right\rangle}

\begin{document}

\title{Sub-Riemannian Distances on Lie Groups: Tangent Asymptotics and Hermitian Matrix Formulas}

\author{Wolfram Bauer, Zhicheng Han, Zhipeng Yang\thanks{Corresponding author: yangzhipeng326@163.com}}
\date{}
\maketitle
\AtEndDocument{%
\par
\bigskip
\bigskip
\noindent
	\textbf{Wolfram Bauer:}\\[0.2em]
	\textsc{Institut f\"ur Analysis, Leibniz Universit\"at Hannover,
30167 Hannover, Germany}\\[0.3em]
	\textit{E-mail address}: \texttt{bauer@math.uni-hannover.de}\\[1.5em]
	\noindent
  \textbf{Zhicheng Han:}\\[0.2em]
	\textsc{Institut f\"ur Analysis, Leibniz Universit\"at Hannover,
30167 Hannover, Germany}\\[0.3em]
	\textit{E-mail address}: \texttt{hanzc@math.uni-hannover.de}\\[1.5em]
  \noindent
	\textbf{Zhipeng Yang:}\\[0.2em]
	\textsc{Department of Mathematics, Yunnan Key Laboratory of Modern Analytical Mathematics and Applications, Yunnan Normal University, 650500 Kunming, China}\\
	\textsc{Department of Mathematics: Analysis, Logic and Discrete Mathematics, Ghent University, 9000 Ghent, Belgium}\\[0.3em]
	\textit{E-mail address}: \texttt{yangzhipeng326@163.com}%
}

\begin{abstract}
Let $d_{\sr}$ be the Carnot--Carath\'{e}odory distance associated with a two-step Lie pair and let $d_0$ be the distance on its nilpotent tangent group. For compact Lie pairs we prove, locally uniformly in $(X,Z)$,
\[
 d_{\sr}\!\left(e,\exp(rX+r^2Z)\right)^2
 =r^2d_0(e,(X,Z))^2+O(r^3).
\]
For symmetric pairs, the remainder improves to $O(r^4)$. This order is sharp in general, and the estimate remains uniform at cut points, abnormal endpoints, and points where the tangent endpoint map changes rank. If the tangent minimizer is unique, the endpoint differential is onto, and the constrained Hessian is coercive, we determine the coefficient of $r^4$ with remainder $O(r^6)$. For block subgroups of $\SU(N)$, we express the vertical squared distance on the tangent group as the attained minimum of a finite-dimensional functional on Hermitian matrices subject to linear constraints on their partial traces. Every minimizer gives a minimizing loop represented by a trigonometric polynomial. For one $m$-dimensional block and $q$ scalar blocks, the coefficient has weights $\lceil j/q\rceil$ and reduces to a trace norm formula when $q\ge\lfloor m/2\rfloor$.
\end{abstract}



\section{Introduction and main results}
\label{Sec1}

The nilpotent tangent group is the basic local model of an equiregular sub-Riemannian manifold. In privileged coordinates, anisotropic dilations rescale the Carnot--Carath\'{e}odory distance to the distance on a Carnot group. Mitchell proved that this group is the pointed metric tangent \cite{Mitchell1985}. See also Bella\"iche \cite{Bellaiche1996} and the monographs \cite{AgrachevBarilariBoscain2020,Jean2014}. This result determines the leading metric geometry near a point. It does not, however, give a uniform estimate for the difference between the original distance and its tangent model.

Such an estimate is delicate at singular endpoints. The distance may fail to be smooth at the cut locus, minimizing controls may be nonunique or abnormal, and the differential of the endpoint map may change rank. An expansion obtained from one minimizing branch need not be uniform near these sets. At a regular endpoint, a further question is to determine the first coefficient beyond the tangent term. These problems are relevant both to the local geometry of horizontal curves and to the logarithmic small-time behavior of hypoelliptic heat kernels.

Quantitative forms of nilpotent approximation were developed by Bella\"iche and Jean \cite{Bellaiche1996,Jean2014}. Le Donne, Nicolussi Golo, and Tettamanti recently obtained a quantitative version of Mitchell's tangent theorem, with rates determined by algebraic parameters of the chosen splitting \cite[Theorem~1.2]{LeDonneNicolussiGoloTettamanti2026}. The lower estimates below also use the vertical perturbation result of Le Donne, Nalon, Nicolussi Golo, and Ryoo. Their estimate implies local Euclidean Lipschitz continuity of the squared distance on every two-step tangent group \cite[Corollary~3.2]{LeDonneNalonNicolussiGoloRyoo2025}.

In \cite{BauerHanYang2026}, we studied heat kernels, vertical distances, and spectral invariants on compact Lie pairs. For compact symmetric pairs we proved
\[
 d_{\sr}(e,\exp(\varepsilon Z))^2
 =\varepsilon\mathfrak F_{G,K}(Z)+O_R(\varepsilon^{3/2}),
\]
together with a lower estimate whose error is $O_R(\varepsilon^2)$. We also expressed $\mathfrak F_{G,K}$ as the attained minimum of a weighted singular-value functional on the free two-step group, used the central quotient to obtain the corresponding formula on the tangent group, and evaluated the one-column AIII family. Those results determine the leading term at vertical endpoints.

The purpose of the present paper is to obtain a quantitative expansion at arbitrary tangent endpoints and to evaluate the vertical squared distance on the tangent group for general block subgroups of $\SU(N)$. For compact Lie pairs, we obtain a uniform cubic remainder for the squared distance. In the symmetric setting, the grading improves this remainder to quartic order. At regular endpoints we determine the coefficient of $r^4$ by a constrained variational argument. For block subgroups, the vertical squared distance is reduced to a finite-dimensional minimization over Hermitian matrices with prescribed partial traces.

Let $G$ be a compact connected Lie group, let $L\subset G$ be a closed connected subgroup, and fix an $\Ad(G)$-invariant inner product on $\gLie$. Write
\[
 \gLie=\lLie\oplus\pLie,
 \qquad \lLie=\Lie(L),
 \qquad \pLie=\lLie^\perp,
\]
and assume
\[
 \gLie=\pLie+[\pLie,\pLie].
\]
The left translates of $\pLie$ define a bracket-generating horizontal distribution, and we write $d_{\sr}$ for its Carnot--Carath\'{e}odory distance. The nilpotent tangent group at the identity is the two-step group $\mathbb N_{G,L}$ on $\pLie\oplus\lLie$. Its second layer is central, and its bracket is determined by
\[
 B:\Lambda^2\pLie\longrightarrow\lLie,
 \qquad B(X\wedge Y)=\pr_{\lLie}[X,Y].
\]
Let $d_0$ denote its Carnot--Carath\'{e}odory distance.

The first result gives a uniform comparison between the two squared distances.

\begin{theorem}
\label{Thm1.1}
For every compact $\mathcal K\subset\pLie\times\lLie$, there are $C_{\mathcal K},r_{\mathcal K}>0$ such that
\begin{equation}\label{eq1.1}
 \left|
 d_{\sr}\!\left(e,\exp(rX+r^2Z)\right)^2
 -r^2d_0(e,(X,Z))^2
 \right|\le C_{\mathcal K}r^3
\end{equation}
for $(X,Z)\in\mathcal K$ and $0<r<r_{\mathcal K}$.

Suppose more generally that $G$ is a connected Lie group and that $\gLie=\kLie\oplus\pLie$ is a symmetric decomposition,
\[
 [\kLie,\kLie]\subset\kLie,
 \qquad [\kLie,\pLie]\subset\pLie,
 \qquad [\pLie,\pLie]\subset\kLie,
\]
whose infinitesimal $\kLie$-action on $\pLie$ is orthogonal. Let $d_{\sr}$ be the left-invariant Carnot--Carath\'{e}odory distance defined by this metric. If $\gLie=\pLie+[\pLie,\pLie]$, let $d_0$ denote the distance on the corresponding tangent group. Then, for every compact $\mathcal K\subset\pLie\times\kLie$, there are $C_{\mathcal K},r_{\mathcal K}>0$ such that
\begin{equation}\label{eq1.2}
 \left|
 d_{\sr}\!\left(e,\exp(rX+r^2Z)\right)^2
 -r^2d_0(e,(X,Z))^2
 \right|\le C_{\mathcal K}r^4
\end{equation}
for $(X,Z)\in\mathcal K$ and $0<r<r_{\mathcal K}$.
\end{theorem}

The upper estimates are obtained by correcting an arbitrary tangent minimizing control with explicit Lie polynomials. The estimates are therefore uniform at cut points, abnormal endpoints, and points where the tangent endpoint map changes rank. The lower estimates follow by rescaling minimizing controls on $G$ and using the local Lipschitz continuity of $d_0^2$. In the symmetric case, the grading separates odd horizontal terms from even vertical terms and gives the additional cancellation in \eqref{eq1.2}.

We compare these orders with \cite[Theorem~1.2]{LeDonneNicolussiGoloTettamanti2026}. Their estimate
\[
 |d(p,q)-d_0(p,q)|
 \le C\max\{d_0(e,p),d_0(e,q)\}^{1+\alpha_0/\beta_{\mathrm{MT}}}
\]
gives an error $O(r^{2+\alpha_0/\beta_{\mathrm{MT}}})$ for the squared distance along $q_r=\exp(rX+r^2Z)$. For a two-step splitting that is not a Carnot grading, $\beta_{\mathrm{MT}}\in\{1,2\}$. If the first nonhomogeneous bracket has a nonzero $\pLie$-component, then $\alpha_0=1$, and the corresponding orders are $O(r^{5/2})$ and $O(r^3)$. For a symmetric grading, $\alpha_0=2$, and the orders are $O(r^3)$ and $O(r^4)$. Thus \eqref{eq1.1} and \eqref{eq1.2} improve the cases $\beta_{\mathrm{MT}}=2$ and have the same order in the cases $\beta_{\mathrm{MT}}=1$.

In the compact symmetric case, let $\theta$ be the involution, put $K=(G^\theta)_0$, and set $L=K$. Taking $X=0$ and $r=\sqrt\varepsilon$ in \eqref{eq1.2} gives, uniformly for $\mathfrak F_{G,K}(Z)\le R$,
\[
 d_{\sr}(e,\exp(\varepsilon Z))^2
 =\varepsilon\mathfrak F_{G,K}(Z)+O_R(\varepsilon^2).
\]
This improves the error $O_R(\varepsilon^{3/2})$ in \cite[Theorem~4.12]{BauerHanYang2026}. The estimate remains uniform at singular vertical endpoints and follows from the same result that treats nonvertical endpoints.

We next determine the coefficient of $r^4$ at regular data. Let $(X,Z)\in\pLie\times\kLie$, and let $u_0\in L^2([0,1];\pLie)$ be the unique tangent minimizing control for $(X,Z)$. Denote the tangent endpoint map by $F_0$, and let $\lambda_0$ be determined by
\[
 DE(u_0)=DF_0(u_0)^*\lambda_0,
 \qquad
 E(u)=\int_0^1|u(s)|^2\,\mathrm ds.
\]
Assume that $DF_0(u_0)$ is onto and that the constrained Hessian is coercive. Write
\[
 F_2(u_0)=\bigl(Q(u_0),R(u_0)\bigr)
\]
for the cubic horizontal and quartic vertical terms in the logarithmic endpoint expansion.

\begin{theorem}
\label{Thm1.2}
Under the preceding uniqueness, regularity, and coercivity assumptions,
\begin{equation}\label{eq1.3}
 d_{\sr}\!\left(e,\exp(rX+r^2Z)\right)^2
 =r^2d_0(e,(X,Z))^2
 -r^4\lambda_0\!\left(Q(u_0),R(u_0)\right)
 +O(r^6).
\end{equation}
For all sufficiently small $r>0$, the rescaled minimizing control is unique and is real analytic as a function of $r^2$.
\end{theorem}

The coefficient in \eqref{eq1.3} is determined by the Lagrange multiplier and the first two terms beyond the tangent endpoint map. The assumptions hold, for example, at the horizontal data considered in \Cref{Prop3.8}. The quartic order in \eqref{eq1.2} is sharp in general. For the standard symmetric pair in $\SU(2)$, the exact vertical formula gives
\[
 d_{\sr}\!\left(e,\exp(r^2zk)\right)^2
 =2\pi|z|r^2-\frac12z^2r^4,
\]
so the coefficient of $r^4$ is nonzero at vertical cut points.

We next obtain an explicit finite-dimensional formula for the vertical squared distance on the tangent group of block subgroups. Put $\kappa=\lfloor\dim\pLie/2\rfloor$. For $\Omega\in\Lambda^2\pLie$, let $\sigma_1(\Omega)\ge\cdots\ge\sigma_\kappa(\Omega)\ge0$ be the nonnegative normal form coefficients of the associated skew-symmetric operator, and set
\[
 \Psi_{\pLie}(\Omega)
 =4\pi\sum_{j=1}^{\kappa}j\,\sigma_j(\Omega).
\]
The results of \cite{BauerHanYang2026} give
\begin{equation}\label{eq1.4}
 \mathfrak F_{G,L}(Z)
 :=d_0(e,(0,Z))^2
 =\min_{B\Omega=Z}\Psi_{\pLie}(\Omega).
\end{equation}
The minimum is attained, $\mathfrak F_{G,L}$ is globally Lipschitz and equivalent to the Euclidean norm on $\lLie$, and every minimizing lift gives a minimizing loop involving at most $\kappa$ positive Fourier frequencies. We use \eqref{eq1.4} to evaluate the affine fibers for block Lie pairs.

Fix a partition
\[
 N=n_1+\cdots+n_s,
\]
let $Z_{\mathbf n}$ be the traceless block scalar torus, choose a closed connected subtorus $T_0\subset Z_{\mathbf n}$, and put
\[
 L_{\mathbf n,T_0}
 =\left(\prod_{a=1}^s\SU(n_a)\right)T_0\subset\SU(N).
\]
Let $\lLie=\Lie(L_{\mathbf n,T_0})$ and let $\mathcal V_{\mathbf n,T_0}^{\mathrm H}$ be the space of block diagonal Hermitian matrices $H$ for which $iH\in\lLie$. For each pair $a<b$, let $S_{ab}\in\Herm(M_{n_a\times n_b}(\C))$. The map $\mathscr Q_{\mathbf n,T_0}$ is obtained by taking the signed block sums of the two partial traces, with the transpose convention in \eqref{eq4.10}, and projecting onto $\mathcal V_{\mathbf n,T_0}^{\mathrm H}$.

For $S\in\Herm(d)$, write $\lambda_1(S)\ge\cdots\ge\lambda_d(S)$ for its eigenvalues. As in \cite{BauerHanYang2026}, let $\lambda_j^+(S)$ and $\lambda_j^-(S)$ denote the positive eigenvalues of $S\in\Herm(d)$ and the absolute values of its negative eigenvalues, respectively, both in decreasing order. If their numbers are $p(S)$ and $q(S)$, set
\[
 \Phi_d(S)
 :=4\pi\left(
 \sum_{j=1}^{p(S)}j\lambda_j^+(S)
 +\sum_{j=1}^{q(S)}j\lambda_j^-(S)
 \right).
\]

\begin{theorem}
\label{Thm1.3}
For every $H\in\mathcal V_{\mathbf n,T_0}^{\mathrm H}$,
\begin{equation}\label{eq1.5}
 d_0(e,(0,iH))^2
 =\min_{\mathscr Q_{\mathbf n,T_0}(S)=H}
 \sum_{a<b}\Phi_{n_an_b}(S_{ab}).
\end{equation}
The minimum is attained. A minimizing family $(S_{ab})$ gives a tangent minimizing loop, and the component corresponding to $a<b$ uses no Fourier frequency larger than $n_an_b$. The right-hand side defines a positively homogeneous and globally Lipschitz function that is equivalent to the Hilbert--Schmidt norm on $\mathcal V_{\mathbf n,T_0}^{\mathrm H}$.
\end{theorem}

Formula \eqref{eq1.5} converts the minimization over horizontal loops into a finite-dimensional problem for Hermitian matrices. The geometry of the block pair appears in the linear constraint, while the function $\Phi_d$ is independent of the partition. The spectral decomposition of a minimizing matrix gives a minimizing loop represented by a trigonometric polynomial. The formula extends the one-column AIII calculation in \cite[Theorem~6.2]{BauerHanYang2026} to arbitrary block families.

For scalar blocks and the full block scalar torus, \eqref{eq1.5} reduces to an $\ell^1$ transportation problem. Consider also
\[
 \mathbf n=(m,\underbrace{1,\ldots,1}_{q\ \mathrm{times}}),
 \qquad
 \widehat H_A=\diag(A,0,\ldots,0),
 \qquad
 A\in\Herm(m),\quad \Tr A=0.
\]
If $r_m=\lfloor m/2\rfloor$, then
\begin{equation}\label{eq1.6}
 \mathfrak F_{\mathbf n,\{e\}}(\widehat H_A)
 =4\pi\sum_{j=1}^{r_m}
 \left\lceil\frac jq\right\rceil
 \bigl(\lambda_j(A)-\lambda_{m+1-j}(A)\bigr).
\end{equation}
When $q\ge r_m$, all weights in \eqref{eq1.6} are equal to one, and
\[
 \mathfrak F_{\mathbf n,\{e\}}(\widehat H_A)
 =4\pi\min_{t\in\R}\Tr|A+tI_m|.
\]
Thus the formula passes from the weighted eigenvalue sum to a trace norm when the number of scalar blocks is sufficiently large.

Prescribed holonomy problems originate in work of Brockett and Montgomery \cite{Brockett1982,Montgomery1990}. Related $K+P$ geometries, control problems on Stiefel and Grassmann manifolds, and holonomic quantum computation are studied in \cite{TanimuraNakaharaHayashi2005,AutenriedMarkina2014,AlbertiniDAlessandro2018,AlbertiniDAlessandroSheller2020,Sonnerborn2024,SonnerbornProjective2025,Sonnerborn2025}. In this setting, \eqref{eq1.5} gives the infinitesimal cost of a prescribed block holonomy, and \eqref{eq1.6} describes the effect of the scalar blocks.

The distance estimates also give logarithmic heat kernel asymptotics. For the horizontal sub-Laplacian $\Delta_{\pLie}=\sum_j\widetilde X_j^{\,2}$, where $\widetilde X_j$ are the left-invariant fields associated with an orthonormal basis of $\pLie$, L\'{e}andre's formula identifies the small-time logarithmic rate with $d_{\sr}^2$. Hence \eqref{eq1.1}--\eqref{eq1.3} give the corresponding expansion when the heat time first tends to zero and the endpoint then approaches the identity. At block vertical endpoints, the leading term is the coefficient in \eqref{eq1.5}.

The paper is organized as follows. \Cref{Sec2} introduces the tangent group and proves the local Lipschitz estimate for the squared tangent distance. \Cref{Sec3} proves the cubic and quartic estimates, determines the coefficient of $r^4$ under the assumptions above, and proves sharpness. \Cref{Sec4} recalls the vertical distance formula from \cite{BauerHanYang2026}, obtains the block matrix formula, evaluates several families, and gives the distance and heat kernel consequences.

\section{Two-step Lie pairs and the tangent distance}
\label{Sec2}

We first describe the tangent group and its endpoint map. We then prove the local Euclidean Lipschitz estimate for the squared tangent distance used in the lower bounds of \Cref{Sec3}.

For the compact comparison theorem and the block applications, $G$ is a compact connected Lie group, $L\subset G$ is a closed connected subgroup, and $\ip{\cdot}{\cdot}$ is an $\Ad(G)$-invariant inner product on $\gLie$. We use the orthogonal splitting
\begin{equation}\label{eq2.1}
 \gLie=\lLie\oplus\pLie,
 \qquad \lLie=\Lie(L),
 \qquad \pLie=\lLie^\perp.
\end{equation}
Invariance and the fact that $\lLie$ is a subalgebra imply
\begin{equation*}
 [\lLie,\pLie]\subset\pLie.
\end{equation*}
The left-invariant horizontal distribution is $\mathcal H_g=(L_g)_*\pLie$. A horizontal curve $\gamma$ based at the identity is described by its left-trivialized control $a\in L^2([0,1];\pLie)$:
\begin{equation*}
 \gamma(0)=e,
 \qquad \gamma(s)^{-1}\dot\gamma(s)=a(s).
\end{equation*}
Its energy is $E(a)=\int_0^1|a(s)|^2\,\mathrm ds$. The square of the Carnot--Carath\'{e}odory distance $d_{\sr}$ is the minimum of this energy among controls on $[0,1]$ with the prescribed endpoint.

We impose
\begin{equation}\label{eq2.2}
 \gLie=\pLie+[\pLie,\pLie].
\end{equation}
Equivalently, the projected bracket
\begin{equation}\label{eq2.3}
 B:\Lambda^2\pLie\longrightarrow\lLie,
 \qquad B(X\wedge Y)=\pr_{\lLie}[X,Y]
\end{equation}
is onto. The nilpotent tangent Lie algebra is $\pLie\oplus\lLie$, with $\lLie$ central and
\begin{equation}\label{eq2.4}
 [(X,Z),(Y,W)]_0=(0,B(X\wedge Y)).
\end{equation}
We identify the simply connected tangent group $\mathbb N_{G,L}$ with $\pLie\times\lLie$ through the exponential map. Its product and dilations are
\begin{align}
 (X,Z)(Y,W)
 &=\left(X+Y,Z+W+\frac12B(X\wedge Y)\right),
 \label{eq2.5}\\
 \delta_r(X,Z)&=(rX,r^2Z).\notag
\end{align}
The first layer carries the restriction of the fixed inner product. We write $d_0$ for the Carnot--Carath\'{e}odory distance on $\mathbb N_{G,L}$. The tangent distance is homogeneous under the dilations in \eqref{eq2.5}: $d_0(e,(rX,r^2Z))^2=r^2d_0(e,(X,Z))^2$.

If $u\in L^2([0,1];\pLie)$ and $x(s)=\int_0^su(\tau)\,\mathrm d\tau$, its tangent endpoint is
\begin{equation}\label{eq2.6}
 X(u)=x(1),
 \qquad
 Z(u)=\frac12\int_0^1B(x(s)\wedge u(s))\,\mathrm ds.
\end{equation}
Consequently,
\begin{equation}\label{eq2.7}
 d_0(e,(X,Z))^2=
 \min\left\{E(u):X(u)=X,\ Z(u)=Z\right\}.
\end{equation}
The minimum is attained because $(\mathbb N_{G,L},d_0)$ is a proper geodesic metric space.

The existence of minimizing controls does not require $G$ to be compact. Choose a complete left-invariant Riemannian metric whose restriction to $\pLie$ is equivalent to the horizontal metric, and denote its distance by $d_{\mathrm R}$. Comparing the Riemannian and horizontal lengths of horizontal curves gives $d_{\mathrm R}\le C d_{\sr}$. The ball--box theorem shows that $d_{\sr}$ induces the manifold topology. See, for example, \cite[Chapter~2]{Jean2014}. Hence a closed $d_{\sr}$-ball is a closed subset of a compact $d_{\mathrm R}$-ball and is therefore compact. Thus $(G,d_{\sr})$ is a proper length space. The metric Hopf--Rinow theorem \cite[Section~2.5]{BuragoBuragoIvanov} gives a minimizing curve between any two points, and constant speed parametrization turns its squared length into the minimum energy.

We use the following consequence of the local ball--box theorem. Fix Euclidean norms on $\pLie$ and $\lLie$.

\begin{lemma}
\label{Lem2.1}
There are a neighborhood $U$ of $0$ in $\gLie$ and a constant $C$ such that
\begin{equation*}
 d_{\sr}(e,\exp(P+Z))
 \le C\bigl(|P|+|Z|^{1/2}\bigr)
\end{equation*}
for $P\in\pLie$, $Z\in\lLie$, and $P+Z\in U$. The analogous estimate holds in $\mathbb N_{G,L}$, globally by homogeneity.
\end{lemma}

\begin{proof}
This is the two-step ball--box theorem in exponential privileged coordinates. See, for example, \cite[Chapter~2]{Jean2014}. The coordinates associated with \eqref{eq2.1} have weights one on $\pLie$ and two on $\lLie$, because \eqref{eq2.3} is onto. The tangent statement follows from the same local estimate and homogeneity.
\end{proof}

The distance $d_0$ need not be differentiable at central or abnormal endpoints. What the lower comparison requires is the weaker, but uniform, Euclidean Lipschitz regularity of its square.

\begin{proposition}
\label{Prop2.2}
For every compact $\mathcal K\subset\pLie\times\lLie$, there is $L_{\mathcal K}>0$ such that
\begin{equation}\label{eq2.8}
 |d_0(e,(X,Z))^2-d_0(e,(X',Z'))^2|
 \le L_{\mathcal K}\bigl(|X-X'|+|Z-Z'|\bigr)
\end{equation}
for all $(X,Z),(X',Z')\in\mathcal K$.
\end{proposition}

\begin{proof}
The key ingredient is the improved vertical perturbation estimate \cite[Corollary~3.2]{LeDonneNalonNicolussiGoloRyoo2025}: there is $C_0>0$ such that
\begin{equation}\label{eq2.9}
 d_0(e,q\,(0,\zeta))^2
 \le d_0(e,q)^2+C_0|\zeta|
\end{equation}
for every $q\in\mathbb N_{G,L}$ and $\zeta\in\lLie$. Applying \eqref{eq2.9} to $q\,(0,\zeta)$ with $-\zeta$ gives the corresponding estimate in the reverse direction.

Put $h=X'-X$. By the triangle inequality and the horizontal segment of length $|h|$,
\begin{equation}\label{eq2.10}
 d_0(e,(X,Z)(h,0))
 \le d_0(e,(X,Z))+|h|.
\end{equation}
On a fixed compact neighborhood of $\mathcal K$, the distances on the right are uniformly bounded. Squaring \eqref{eq2.10}, and then reversing the roles of the endpoints, shows that the squared distance changes by at most $C_{\mathcal K}|h|$. Finally,
\[
 (X',Z')=(X,Z)(h,0)(0,\zeta),
 \qquad
 \zeta=Z'-Z-\frac12B(X\wedge h).
\]
On $\mathcal K$, $|\zeta|\le C_{\mathcal K}(|Z'-Z|+|h|)$. Combining this observation with \eqref{eq2.9} and repeating the argument in the reverse direction proves \eqref{eq2.8}.
\end{proof}

\begin{remark}
The squared distance in \Cref{Prop2.2} is locally Euclidean Lipschitz even at the identity. For a fixed central vector $Z\ne0$, homogeneity gives
\[
 d_0(e,(0,\varepsilon Z))
 =\varepsilon^{1/2}d_0(e,(0,Z)),\qquad \varepsilon>0.
\]
Thus the distance has order $\varepsilon^{1/2}$ along this ray, whereas its square has order $\varepsilon$. Away from the identity, \Cref{Prop2.2} also implies local Euclidean Lipschitz continuity of $q\mapsto d_0(e,q)$, since this function is locally bounded away from zero. Local Lipschitz continuity does not imply semiconcavity. For examples involving abnormal minimizers, see \cite{MontanariMorbidelli2016}.
\end{remark}

\section{Quantitative tangent expansions}
\label{Sec3}

We now compare the tangent distance with the distance on the original Lie group. We first establish uniform Magnus expansions and then construct the endpoint corrections used in the upper estimates. The regular case is treated by the implicit function theorem. The final $\SU(2)$ example proves that the quartic order is sharp in general.

\subsection{Uniform endpoint expansions}

We first give the endpoint expansion used throughout the comparison. If $a\in L^1([0,1];\gLie)$, let $\Gamma^a$ solve $(\Gamma^a)^{-1}\dot\Gamma^a=a$, $\Gamma^a(0)=e$. For $\norm{a}_{L^1}$ small, put $\Omega^a=\log\Gamma^a$. The logarithmic differential equation is
\begin{equation}\label{eq3.1}
 \dot\Omega^a
 =a+\frac12[\Omega^a,a]
  +\frac1{12}[\Omega^a,[\Omega^a,a]]+\cdots.
\end{equation}
\begin{lemma}
\label{Lem3.1}
For every $M>0$ and integer $m\ge1$, there are $r_0,C>0$ such that
\begin{equation}\label{eq3.2}
 \left\|
 \log\Gamma^{rb}(1)-\sum_{j=1}^{m}r^j\Omega_j(b)
 \right\|\le C|r|^{m+1}
\end{equation}
whenever $b\in L^1([0,1];\gLie)$, $\|b\|_{L^1}\le M$, and $|r|\le r_0$. Here $\Omega_j$ is the Magnus polynomial of degree $j$. The logarithmic endpoint map is real analytic on a sufficiently small ball in $L^1$. More precisely, for every fixed $\ell\ge0$, after decreasing $r_0$ if necessary,
\begin{equation}\label{eq3.3}
 \sup_{\norm{b}_{L^1}\le M}
 \norm{D_b^\ell\!\left(
  \log\Gamma^{rb}(1)-\sum_{j=1}^{m}r^j\Omega_j(b)
 \right)}
 \le C_{M,m,\ell}|r|^{m+1},
\end{equation}
where the norm is the operator norm of an $\ell$-linear map from $(L^1)^\ell$ to $\gLie$.
\end{lemma}

\begin{proof}
Work in a Baker--Campbell--Hausdorff (BCH) neighborhood of the origin in $\gLie$. By Ado's theorem, there is an injective finite-dimensional representation
\[
 \rho:\gLie\longrightarrow\mathfrak{gl}(N,\R).
\]
After shrinking the neighborhood, $\rho$ integrates to a local Lie group homomorphism $\varrho$ that is an embedding onto its image and satisfies
\[
 \varrho(\exp A)=\exp(\rho(A)).
\]
Choose a smaller exponential neighborhood $U$ on which both the group logarithm and the matrix logarithm are defined. Standard ODE estimates give $\Gamma^a([0,1])\subset U$ whenever $\|a\|_{L^1}\le\delta$, for some $\delta>0$.

Let $Y_a(t)=\varrho(\Gamma^a(t))$. Since $\dot Y_a=Y_a\rho(a)$, the chronological series is
\[
 Y_a(1)=I+\sum_{n=1}^{\infty}P_n(a),
\]
where
\[
 P_n(a)=
 \int_{0<t_1<\cdots<t_n<1}
 \rho(a(t_1))\cdots\rho(a(t_n))
 \,\mathrm dt_1\cdots\,\mathrm dt_n.
\]
For a suitable constant $C_0$, symmetry of the scalar majorant gives
\[
 \|P_n(a)\|
 \le \frac{(C_0\|a\|_{L^1})^n}{n!}.
\]
More generally, for $0\le\ell\le n$, differentiation under the integral and summation over the positions occupied by $h_1,\ldots,h_\ell$ give
\begin{equation}\label{eq3.4}
 \bigl\|D^\ell P_n(a)[h_1,\ldots,h_\ell]\bigr\|
 \le
 \frac{C_0^n\|a\|_{L^1}^{n-\ell}}{(n-\ell)!}
 \prod_{j=1}^{\ell}\|h_j\|_{L^1}.
\end{equation}
To see the factorial, label the $n-\ell$ copies of $a$. The ordered simplices obtained from all positions of the labelled functions partition the full product domain, while permutations of the copies of $a$ give the factor $(n-\ell)!$. Estimate \eqref{eq3.4} shows that the chronological series and each of its Fr\'{e}chet derivatives converge uniformly on bounded subsets of $L^1$. After complexification, it defines an entire analytic map from $L^1([0,1];\gLie_\C)$ to the matrix algebra, where $\gLie_\C=\gLie\otimes_\R\C$.

On a smaller ball, $\|Y_a(1)-I\|<1/2$, so the matrix logarithm is analytic and its power series converges uniformly together with all derivatives. Since
\[
 \rho\!\left(\log\Gamma^a(1)\right)=\log Y_a(1),
\]
injectivity of $\rho$ and equivalence of finite-dimensional norms show that the logarithmic endpoint map
\[
 \Phi(a):=\log\Gamma^a(1)
\]
is analytic on an $L^1$-ball. Its Taylor expansion at the origin has the form
\[
 \Phi(a)=\sum_{j=1}^{\infty}L_j(a,\ldots,a)
\]
where $L_j$ is a continuous symmetric $j$-linear map. Substitution into \eqref{eq3.1}, or equivalently the homogeneous Chen--Strichartz expansion \cite{Strichartz1987}, identifies
\[
 L_j(b,\ldots,b)=\Omega_j(b).
\]

Fix $M,m,\ell$. Choose $R>0$ so that $RM$ lies strictly inside the analytic ball of $\Phi$. The estimates above show that, for $\|b\|_{L^1}\le M$, the map
\[
 z\longmapsto D^\ell\bigl[\Phi(z\,\cdot)\bigr](b)
\]
is holomorphic and uniformly bounded on $|z|\le R$. Its Taylor series is
\[
 D^\ell\bigl[\Phi(z\,\cdot)\bigr](b)
 =\sum_{j\ge1}z^jD^\ell\Omega_j(b).
\]
Cauchy's integral formula, applied on the circle $|z|=R$, gives a constant $C_{M,m,\ell}$ such that, for $|r|\le R/2$,
\[
 \sup_{\|b\|_{L^1}\le M}
 \left\|D_b^\ell\!\left(
 \Phi(rb)-\sum_{j=1}^{m}r^j\Omega_j(b)
 \right)\right\|
 \le C_{M,m,\ell}|r|^{m+1}.
\]
This proves \eqref{eq3.3}. The case $\ell=0$ gives \eqref{eq3.2}.
\end{proof}

We use the Magnus expansion up to order four. Controls bounded in $L^2([0,1])$ are bounded in $L^1([0,1])$. Under the symmetric grading, the horizontal Magnus series contains odd degrees and the vertical series contains even degrees. Hence the first omitted terms in the horizontal and vertical components have degrees five and six. After division by $r$ and $r^2$ in the rescaled endpoint map, both remainders are $O(r^4)$, with the same estimate after any fixed number of control Fr\'{e}chet derivatives.

\subsection{Cubic comparison for compact Lie pairs}

\begin{theorem}
\label{Thm3.2}
For a compact Lie pair satisfying \eqref{eq2.1}--\eqref{eq2.2}, and for every compact set $\mathcal K\subset\pLie\times\lLie$, there are $C_{\mathcal K},r_{\mathcal K}>0$ such that
\begin{equation}\label{eq3.5}
 \left|
 d_{\sr}\!\left(e,\exp(rX+r^2Z)\right)^2-r^2d_0(e,(X,Z))^2
 \right|\le C_{\mathcal K}r^3
\end{equation}
whenever $(X,Z)\in\mathcal K$ and $0<r<r_{\mathcal K}$.
\end{theorem}

\begin{proof}
For the upper bound, fix $(X,Z)\in\mathcal K$ and choose an arbitrary minimizing tangent control $u$ in \eqref{eq2.7}. Set
\begin{equation*}
 x(s)=\int_0^su(\tau)\,\mathrm d\tau,
 \qquad
 P=\frac12\int_0^1\pr_{\pLie}[x,u]\,\mathrm ds,
 \qquad
 h=\pr_{\pLie}[x,u],
 \qquad
 H(s)=\int_0^sh(\tau)\,\mathrm d\tau.
\end{equation*}
Then $x(1)=X$, $H(1)=2P$, and
\[
 \frac12\int_0^1[x,u]\,\mathrm ds=P+Z.
\]
Define
\begin{equation*}
 v=-\frac13h,
 \qquad a_r=r(u+rv).
\end{equation*}
Writing
\[
 \log\Gamma^{a_r}(1)=rx(1)+r^2y(1)+r^3q(1)+O_{\mathcal K}(r^4),
\]
the first two equations obtained from \eqref{eq3.1} are
\[
 \dot x=u,
 \qquad
 \dot y=v+\frac12[x,u].
\]
Consequently,
\begin{equation}\label{eq3.6}
 y_{\pLie}(1)=\frac16H(1)=\frac13P,
 \qquad
 y_{\lLie}(1)=Z.
\end{equation}
At degree three, the Magnus expansion gives
\[
 \dot q=\frac12\bigl([x,v]+[y,u]\bigr)
  +\frac1{12}[x,[x,u]].
\]
Only brackets of two $\pLie$-components contribute to the $\lLie$-projection. Since $v=-h/3$ and $y_{\pLie}=H/6$,
\[
 \pr_{\lLie}\dot q
 =-\frac1{12}\pr_{\lLie}[x,h]
  +\frac1{12}\pr_{\lLie}[H,u].
\]
The identity
\[
 \frac{\mathrm d}{\mathrm ds}[x,H]=[x,h]-[H,u]
\]
therefore gives
\begin{equation}\label{eq3.7}
 \pr_{\lLie}q(1)
 =-\frac1{12}\pr_{\lLie}[X,H(1)]
 =-\frac16\pr_{\lLie}[X,P].
\end{equation}

Place the horizontal correction of order $r^2$ before the main curve and set
\begin{equation*}
 g_r=\exp\!\left(-\frac{r^2}{3}P\right)\Gamma^{a_r}(1).
\end{equation*}
At degree three, the preparatory segment contributes the BCH term
\[
 \frac12\left[-\frac{P}{3},X\right]=\frac16[X,P].
\]
Its $\lLie$-projection cancels \eqref{eq3.7}, while its term of order $r^2$ cancels the horizontal component in \eqref{eq3.6}. Hence
\begin{equation*}
 \log g_r=rX+r^2Z+E_r^{\pLie}+E_r^{\lLie},
 \qquad
 E_r^{\pLie}=O_{\mathcal K}(r^3),\quad E_r^{\lLie}=O_{\mathcal K}(r^4).
\end{equation*}
Let $T_r=rX+r^2Z$ and $E_r=E_r^{\pLie}+E_r^{\lLie}$. The BCH series gives
\[
 \log\!\left(g_r^{-1}\exp(T_r)\right)
 =\operatorname{BCH}(-T_r-E_r,T_r)
 =-E_r-\frac12[E_r,T_r]+O_{\mathcal K}(r^5).
\]
Since $[\lLie,\pLie]\subset\pLie$, the horizontal component of this logarithm is $O_{\mathcal K}(r^3)$ and its vertical component is $O_{\mathcal K}(r^4)$. By \Cref{Lem2.1}, the relative element can be joined to the identity by a horizontal curve of length $O_{\mathcal K}(r^2)$. Left translation of this closing curve from $g_r$ completes the concatenation at $\exp(rX+r^2Z)$.

The constants are uniform without a continuous choice of tangent minimizer. Indeed,
\begin{equation*}
 E(u)=d_0(e,(X,Z))^2,\qquad
 \norm{x}_{L^\infty}\le\norm{u}_{L^2},\qquad
 |P|+\norm{v}_{L^2}+\norm{y}_{L^\infty}+|q(1)|\le C_{\mathcal K},
\end{equation*}
because the Magnus coefficients are finite Lie polynomials in $u$ and the squared tangent distance is bounded on $\mathcal K$. Invariance of the inner product gives
\[
 \ip{u}{v}=-\frac13\ip{u}{[x,u]}=0
\]
pointwise. Therefore
\begin{equation*}
 E(a_r)=r^2d_0(e,(X,Z))^2+r^4E(v).
\end{equation*}
The main curve has length at most $r d_0(e,(X,Z))+r^2\sqrt{E(v)}$. The preparatory and closing curves have length $O_{\mathcal K}(r^2)$. Squaring the resulting length bound, using the boundedness of the squared tangent distance on $\mathcal K$, gives
\[
 d_{\sr}\!\left(e,\exp(rX+r^2Z)\right)^2
 \le r^2d_0(e,(X,Z))^2+C_{\mathcal K}r^3.
\]
This estimate also covers $d_0(e,(X,Z))^2=0$, since no division by $d_0(e,(X,Z))$ occurs.

For the lower bound, let $a_r$ be a minimizing control on $[0,1]$ for the endpoint $\exp(rX+r^2Z)$, and put $u_r=a_r/r$. The upper bound makes $E(u_r)$ uniformly bounded. Define
\[
 x_r(s)=\int_0^su_r(\tau)\,\mathrm d\tau,
 \qquad
 P_r+Z_r=\frac12\int_0^1[x_r,u_r]\,\mathrm ds,
 \quad P_r\in\pLie,\quad Z_r\in\lLie.
\]
The expansion through degree two is
\[
 \log\Gamma^{a_r}(1)
 =r x_r(1)+r^2(P_r+Z_r)+O_{\mathcal K}(r^3).
\]
Comparison with the exact logarithm $rX+r^2Z$, together with the uniform bounds on $P_r$, gives
\begin{equation}\label{eq3.8}
 x_r(1)=X+O_{\mathcal K}(r),\qquad Z_r=Z+O_{\mathcal K}(r).
\end{equation}
The control $u_r$ is a tangent competitor for $(x_r(1),Z_r)$, so
\[
 d_0(e,(x_r(1),Z_r))^2\le E(u_r)
 =r^{-2}d_{\sr}\!\left(e,\exp(rX+r^2Z)\right)^2.
\]
After enlarging $\mathcal K$ slightly, \Cref{Prop2.2} and \eqref{eq3.8} yield
\[
 d_0(e,(X,Z))^2\le E(u_r)+C_{\mathcal K}r.
\]
Multiplication by $r^2$ proves the lower bound in \eqref{eq3.5}.
\end{proof}

The homogeneous version of the result is sometimes more convenient. Put

\[
 \rho(P+Z)=|P|+|Z|^{1/2},\qquad P\in\pLie,\quad Z\in\lLie.
\]

\begin{corollary}
\label{Cor3.3}
There are a neighborhood $U$ of $0$ in $\gLie$ and $C>0$ such that
\begin{equation*}
 \left|d_{\sr}(e,\exp(P+Z))^2-d_0(e,(P,Z))^2\right|
 \le C\rho(P+Z)^3
\end{equation*}
for $P+Z\in U$.
\end{corollary}

\begin{proof}
For $\rho>0$, apply \Cref{Thm3.2} with $r=\rho$, $X=P/r$, and $W=Z/r^2$. These normalized points belong to the fixed compact set $\{|X|\le1,|W|\le1\}$. Use the homogeneity of the squared tangent distance. The case $\rho=0$ is immediate.
\end{proof}

\subsection{Quartic comparison in the symmetric setting}
\label{Subsec3.3}

The quartic argument requires less global structure than the preceding compact comparison theorem. We assume that the Lie algebra has a symmetric decomposition
\begin{equation}\label{eq3.9}
 \gLie=\kLie\oplus\pLie,\qquad
 [\kLie,\kLie]\subset\kLie,\qquad
 [\kLie,\pLie]\subset\pLie,\qquad [\pLie,\pLie]\subset\kLie,
\end{equation}
and that the inner product on $\pLie$ is invariant under the infinitesimal $\kLie$-action:
\begin{equation}\label{eq3.10}
 \ip{[A,X]}{Y}+\ip{X}{[A,Y]}=0
 \qquad(A\in\kLie,\ X,Y\in\pLie).
\end{equation}
We assume $\gLie=\pLie+[\pLie,\pLie]$. The projected bracket is then the bracket itself, and it defines the tangent group $\mathbb N_{G,K}$ and its distance $d_0$ exactly as in \eqref{eq2.4}--\eqref{eq2.7}, with $\lLie$ replaced by $\kLie$. Every symmetric pair in the compact setting of \Cref{Sec2} has this structure. Only the Lie algebraic assumptions in \eqref{eq3.9}--\eqref{eq3.10} are used, so the argument applies to connected compact and noncompact groups.

\begin{theorem}
\label{Thm3.4}
Let $G$ be a connected Lie group whose Lie algebra and horizontal inner product satisfy \eqref{eq3.9}--\eqref{eq3.10}, and assume $\gLie=\pLie+[\pLie,\pLie]$. For every compact $\mathcal K\subset\pLie\times\kLie$, there are $C_{\mathcal K},r_{\mathcal K}>0$ such that
\begin{equation}\label{eq3.11}
 \left|
 d_{\sr}\!\left(e,\exp(rX+r^2Z)\right)^2-r^2d_0(e,(X,Z))^2
 \right|\le C_{\mathcal K}r^4
\end{equation}
for $(X,Z)\in\mathcal K$ and $0<r<r_{\mathcal K}$.
\end{theorem}

The norm equivalence in \Cref{Prop4.1} makes every sublevel of $\mathfrak F_{G,K}$ compact. Hence, for $X=0$ and $r=\sqrt\varepsilon$, \Cref{Thm3.4} gives, uniformly for $\mathfrak F_{G,K}(Z)\le R$,
\[
 d_{\sr}(e,\exp(\varepsilon Z))^2
 =\varepsilon\mathfrak F_{G,K}(Z)+O_R(\varepsilon^2).
\]
This improves the error in absolute value from $O_R(\varepsilon^{3/2})$ in \cite[Theorem~4.12]{BauerHanYang2026} to $O_R(\varepsilon^2)$. The same theorem already gives the lower estimate with error $O_R(\varepsilon^2)$.

\begin{proof}
For the upper bound, choose a tangent minimizing control $u$ for $(X,Z)$, and set
\[
 x(s)=\int_0^su(\tau)\,\mathrm d\tau,
 \qquad
 \dot z=\frac12[x,u],\qquad z(0)=0.
\]
Thus $x(1)=X$ and $z(1)=Z$. For the uncorrected control $ru$, write
\begin{equation*}
 \log\Gamma^{ru}(1)
 =rX+r^2Z+r^3Q+r^4R
  +O_{\mathcal K}(r^5)_{\pLie}+O_{\mathcal K}(r^6)_{\kLie}.
\end{equation*}
The grading places $Q$ in $\pLie$ and $R$ in $\kLie$. The construction below cancels the cubic horizontal coefficient and the quartic vertical coefficient. Parity then leaves an error of order $r^5$ in $\pLie$ and order $r^6$ in $\kLie$. By \Cref{Lem2.1}, this error can be corrected by a horizontal curve of length $O(r^3)$. If $q,w$ are the corresponding homogeneous paths, then
\begin{align}
 \dot q&=\frac12[z,u]+\frac1{12}[x,[x,u]],
 \label{eq3.12}\\
 \dot w&=\frac12[q,u]
 +\frac1{12}\bigl([x,[z,u]]+[z,[x,u]]\bigr),
 \label{eq3.13}
\end{align}
with $Q=q(1)$ and $R=w(1)$.

Introduce
\begin{equation*}
 A=\int_0^1[z,u]\,\mathrm ds,
 \qquad
 \mathcal B=\int_0^1[x,[x,u]]\,\mathrm ds,
 \qquad
 I=\int_0^1[x,[z,u]]\,\mathrm ds.
\end{equation*}
Equation \eqref{eq3.12} gives
\begin{equation}\label{eq3.14}
 Q=\frac12A+\frac1{12}\mathcal B.
\end{equation}
To obtain the identity of degree four, set
\[
 J=\int_0^1[z,[x,u]]\,\mathrm ds,
 \qquad
 K_4=\int_0^1[x,[x,[x,u]]]\,\mathrm ds.
\]
Since
\[
 \frac{\mathrm d}{\mathrm ds}[q,x]=[\dot q,x]+[q,u],
\]
integration and \eqref{eq3.12} give
\[
 \frac12\int_0^1[q,u]\,\mathrm ds
 =-\frac12[X,Q]+\frac14I+\frac1{24}K_4.
\]
Substitution into \eqref{eq3.13} yields
\begin{equation}\label{eq3.15}
 R=-\frac12[X,Q]+\frac13I+\frac1{12}J+\frac1{24}K_4.
\end{equation}
Jacobi's identity gives
\begin{equation}\label{eq3.16}
 [x,[z,u]]=[[x,z],u]+[z,[x,u]].
\end{equation}
Let $L=\int_0^1[[x,z],u]\,\mathrm ds$. Differentiating $[[x,z],x]$ and using $\dot z=[x,u]/2$, we obtain
\begin{align*}
 L
 &= [[X,Z],X]
  -\int_0^1[[u,z],x]\,\mathrm ds
  -\int_0^1[[x,\dot z],x]\,\mathrm ds\\
 &= [X,[Z,X]]-I+\frac12K_4.
\end{align*}
Here $[[u,z],x]=[x,[z,u]]$ and $[[x,\dot z],x]=-[x,[x,[x,u]]]/2$. Integrating \eqref{eq3.16} gives $I=L+J$, and hence
\begin{equation}\label{eq3.17}
 J=2I-\frac12K_4-[X,[Z,X]].
\end{equation}
Combining \eqref{eq3.15} and \eqref{eq3.17} cancels the $K_4$-terms and gives
\begin{equation}\label{eq3.18}
 R=\frac12I-\frac12[X,Q]-\frac1{12}[X,[Z,X]].
\end{equation}

Choose the control correction
\begin{equation*}
 \eta=-\frac12[z,u],
 \qquad
 a_r=r(u+r^2\eta).
\end{equation*}
For the tangent endpoint map, its first variations are
\begin{align}
 \delta X&=\int_0^1\eta\,\mathrm ds=-\frac12A,
 \label{eq3.19}\\
 \delta Z&=\int_0^1[x,\eta]\,\mathrm ds-\frac12[X,\delta X]
  =-\frac12I+\frac14[X,A].
 \label{eq3.20}
\end{align}
The first equality in \eqref{eq3.20} follows by integrating $[\delta x,u]$ by parts. Equations \eqref{eq3.14}, \eqref{eq3.18}, and \eqref{eq3.19}--\eqref{eq3.20} show that the main curve has residual coefficients
\begin{equation*}
 H_3=\frac1{12}\mathcal B,
 \qquad
 V_4=-\frac1{24}[X,\mathcal B]-\frac1{12}[X,[Z,X]].
\end{equation*}

Define
\begin{equation*}
 C=-\frac1{12}\bigl(\mathcal B+[Z,X]\bigr),
 \qquad
 D=\frac1{12}[Z,X],
\end{equation*}
and concatenate
\begin{equation*}
 g_r=\exp(r^3C)\,\Gamma^{a_r}(1)\,\exp(r^3D).
\end{equation*}
Both $C$ and $D$ lie in $\pLie$. Since $C+D=-H_3$, the cubic horizontal residual vanishes. At degree four, the two segments contribute to the vertical component the term
\[
 \frac12[X,D-C]
 =\frac1{24}[X,\mathcal B]+\frac1{12}[X,[Z,X]]=-V_4.
\]
Thus
\begin{equation*}
 \log g_r=rX+r^2Z+O_{\mathcal K}(r^5)_{\pLie}+O_{\mathcal K}(r^5)_{\kLie}.
\end{equation*}

The remaining vertical order is determined by the grading. The map
\[
 \theta(A+Y)=A-Y,
 \qquad A\in\kLie,
 \quad Y\in\pLie,
\]
is a Lie algebra involution. It integrates on a BCH neighborhood to a local Lie group automorphism $\Theta$ satisfying $\Theta(\exp W)=\exp(\theta W)$. Since $\eta,C,D\in\pLie$,
\[
 \theta(a_r)=a_{-r},\qquad
 \Theta(\exp(r^3C))=\exp((-r)^3C),\qquad
 \Theta(\exp(r^3D))=\exp((-r)^3D).
\]
Uniqueness of solutions of the control equation gives $\Theta(\Gamma^{a_r})=\Gamma^{a_{-r}}$, and therefore
\begin{equation}\label{eq3.21}
 \Theta(g_r)=g_{-r},
 \qquad
 \theta(\log g_r)=\log g_{-r}.
\end{equation}
Writing $\log g_r=P(r)+L(r)$ with $P(r)\in\pLie$, $L(r)\in\kLie$, equation \eqref{eq3.21} implies
\[
 P(-r)=-P(r),\qquad L(-r)=L(r).
\]
The vertical component has no term of degree five. Hence
\begin{equation}\label{eq3.22}
 \log g_r=rX+r^2Z+O_{\mathcal K}(r^5)_{\pLie}+O_{\mathcal K}(r^6)_{\kLie}.
\end{equation}
This proves the required parity relation.

Applying the BCH formula to $g_r^{-1}\exp(rX+r^2Z)$ preserves the orders in \eqref{eq3.22}: its horizontal logarithmic component is $O_{\mathcal K}(r^5)$ and its vertical component is $O_{\mathcal K}(r^6)$. By \Cref{Lem2.1}, a logarithmic error with horizontal component $O_{\mathcal K}(r^5)$ and vertical component $O_{\mathcal K}(r^6)$ can be closed by a horizontal curve of length $O_{\mathcal K}(r^3)$.

All correction terms are uniformly bounded on $\mathcal K$. Indeed,
\[
 \norm{z}_{L^\infty}\le C\norm{u}_{L^2}^2,
 \qquad
 \norm{\eta}_{L^2}+|A|+|\mathcal B|+|I|+|C|+|D|\le C_{\mathcal K}.
\]
The infinitesimal orthogonality in \eqref{eq3.10} gives $\ip{u}{\eta}=0$ pointwise, and hence
\begin{equation*}
 E(a_r)=r^2d_0(e,(X,Z))^2+r^6E(\eta).
\end{equation*}
The main curve has length at most $r d_0(e,(X,Z))+r^3\sqrt{E(\eta)}$. The two horizontal segments and the closing curve have length $O_{\mathcal K}(r^3)$. After squaring, this proves the upper bound in \eqref{eq3.11}.

For the lower bound, let $a_r$ be a minimizing control on $G$, put $u_r=a_r/r$, and use the upper bound to obtain $\sup_r\|u_r\|_{L^2}<\infty$. Write $(X_r,Z_r)$ for the tangent endpoint of $u_r$. Uniformly on this bounded set, the Magnus expansion and the parity relation give
\begin{align*}
 \pr_{\pLie}\log\Gamma^{ru_r}(1)
 &=rX_r+r^3Q(u_r)+O_{\mathcal K}(r^5),\\
 \pr_{\kLie}\log\Gamma^{ru_r}(1)
 &=r^2Z_r+r^4R(u_r)+O_{\mathcal K}(r^6).
\end{align*}
The Lie polynomials $Q(u_r)$ and $R(u_r)$ are uniformly bounded. Comparison with the exact logarithm $rX+r^2Z$ yields
\begin{equation}\label{eq3.23}
 X_r=X-r^2Q(u_r)+O_{\mathcal K}(r^4)=X+O_{\mathcal K}(r^2),
 \qquad
 Z_r=Z-r^2R(u_r)+O_{\mathcal K}(r^4)=Z+O_{\mathcal K}(r^2).
\end{equation}
Since $u_r$ is a tangent competitor for $(X_r,Z_r)$, $d_0(e,(X_r,Z_r))^2\le E(u_r)$. \Cref{Prop2.2} and \eqref{eq3.23} give
\[
 d_0(e,(X,Z))^2\le E(u_r)+C_{\mathcal K}r^2.
\]
Multiplication by $r^2$ proves the lower bound with remainder $C_{\mathcal K}r^4$.
\end{proof}

\begin{corollary}
\label{Cor3.5}
Under the hypotheses of \Cref{Thm3.4}, there are a neighborhood $U$ of $0$ in $\gLie$ and $C>0$ such that
\begin{equation*}
 \left|d_{\sr}(e,\exp(P+Z))^2-d_0(e,(P,Z))^2\right|
 \le C\rho(P+Z)^4.
\end{equation*}
\end{corollary}

\begin{proof}
Normalize by $r=\rho(P+Z)$ exactly as in \Cref{Cor3.3} and apply \Cref{Thm3.4}.
\end{proof}

\subsection{The coefficient of order four at regular data}

The estimate in \Cref{Thm3.4} is uniform across the tangent cut locus. We now assume that the tangent minimizer is unique, the endpoint constraint is regular, and the constrained second variation is coercive. Under these assumptions the coefficient of $r^4$ is well defined. Put
\[
 \mathcal H=L^2([0,1];\pLie),\qquad \mathcal Y=\pLie\oplus\kLie,
\]
and let $F_0:\mathcal H\to\mathcal Y$ be the tangent endpoint map in \eqref{eq2.6}. For $r\ne0$, define
\begin{equation}\label{eq3.24}
 F(r^2,u)=\left(
 r^{-1}\pr_{\pLie}\log\Gamma^{ru}(1),
 r^{-2}\pr_{\kLie}\log\Gamma^{ru}(1)
 \right).
\end{equation}
The grading involution makes the horizontal logarithmic component odd in $r$ and the vertical component even. Hence the apparent singularities in \eqref{eq3.24} are removable. With $\tau=r^2$, the map $F(\tau,u)$ extends to a real analytic map between Banach spaces near $(0,u_0)$ and satisfies
\[
 F(\tau,u)=F_0(u)+\tau F_2(u)+O(\tau^2),
 \qquad F_2(u)=(Q(u),R(u)).
\]
More precisely, for every $M>0$ and $\ell=0,1,2$, there are $\tau_M,C_{M,\ell}>0$ such that
\begin{equation}\label{eq3.25}
 \sup_{\norm{u}_{\mathcal H}\le M}
 \norm{D_u^\ell\!\left(
 F(\tau,u)-F_0(u)-\tau F_2(u)
 \right)}
 \le C_{M,\ell}|\tau|^2
\end{equation}
whenever $|\tau|\le\tau_M$. This follows from \Cref{Lem3.1} and the continuous embedding $L^2([0,1])\hookrightarrow L^1([0,1])$. Explicitly, if
\[
 x(t)=\int_0^t u(v)\,\mathrm dv,
 \qquad
 \dot z=\frac12[x,u],\qquad z(0)=0,
\]
and $X=x(1)$, $Z=z(1)$, then
\begin{align*}
 Q(u)&=\frac12\int_0^1[z,u]\,\mathrm dt
  +\frac1{12}\int_0^1[x,[x,u]]\,\mathrm dt,
 \\
 R(u)&=\frac12\int_0^1[x,[z,u]]\,\mathrm dt
  -\frac12[X,Q(u)]-\frac1{12}[X,[Z,X]].
\end{align*}

Suppose that $u_0$ minimizes $E$ under $F_0(u)=y$ and that $D_0:=DF_0(u_0)$ is onto. Since $u_0$ is a constrained minimizer, $DE(u_0)[h]=0$ for every $h\in\ker D_0$. Hence
\[
 \lambda_0(D_0h)=DE(u_0)[h]
\]
defines a unique bounded functional $\lambda_0\in\mathcal Y^*$, and
\[
 DE(u_0)=D_0^*\lambda_0.
\]

\begin{assumption}\label{Ass3.6}
We impose the following assumptions at $y=(X,Z)\in\mathcal Y$.
\begin{enumerate}
\item The problem
 \[
  \min\{E(u):F_0(u)=y\}
 \]
has a unique global minimizer $u_0\in\mathcal H$.
\item The differential $D_0=DF_0(u_0):\mathcal H\to\mathcal Y$ is onto.
\item Let $\lambda_0\in\mathcal Y^*$ be determined by
 \begin{equation}\label{eq3.26}
  2\ip{u_0}{h}_{\mathcal H}=\lambda_0(D_0h)
  \qquad(h\in\mathcal H).
 \end{equation}
The constrained Hessian
 \begin{equation}\label{eq3.27}
 H_0(h,k)=2\ip{h}{k}_{\mathcal H}
 -\lambda_0\!\left(D^2F_0(u_0)[h,k]\right)
 \end{equation}
is coercive on $\ker D_0$: there is $c_0>0$ such that
 \[
  H_0(h,h)\ge c_0\norm{h}_{\mathcal H}^2
  \qquad(h\in\ker D_0).
 \]
\end{enumerate}
The last condition is stronger than nondegeneracy of the constrained Hessian.
\end{assumption}

\begin{theorem}\label{Thm3.7}
Let $G$ satisfy the symmetric hypotheses of \Cref{Thm3.4}, and suppose that the conditions in \Cref{Ass3.6} hold at $(X,Z)$. Then
\begin{equation}\label{eq3.28}
 d_{\sr}\!\left(e,\exp(rX+r^2Z)\right)^2
 =r^2d_0(e,(X,Z))^2+r^4c_4(X,Z)+O(r^6),
\end{equation}
where
\begin{equation}\label{eq3.29}
 c_4(X,Z)=-\lambda_0\!\left(Q(u_0),R(u_0)\right).
\end{equation}
For all sufficiently small $r>0$, the rescaled minimizing control is unique and is real analytic as a function of $r^2$.
\end{theorem}

\begin{proof}
Fix $y=(X,Z)$ and use $\tau=r^2$ as the perturbation parameter. Consider the Lagrange multiplier map
\[
 \mathcal M(\tau,u,\lambda)=
 \left(DE(u)-D_uF(\tau,u)^*\lambda,
  F(\tau,u)-y\right).
\]
It is real analytic near $(0,u_0,\lambda_0)$. Its derivative in $(u,\lambda)$ at that point is
\begin{equation}\label{eq3.30}
 (h,\mu)\longmapsto
 \left(A_0h-D_0^*\mu,D_0h\right),
 \qquad
 \ip{A_0h}{k}_{\mathcal H}=H_0(h,k).
\end{equation}
We verify that this operator has a bounded inverse. Since $D_0$ is onto and $\mathcal Y$ is finite-dimensional, choose a bounded right inverse $\mathcal R_0:\mathcal Y\to\mathcal H$. Given $(f,\xi)\in\mathcal H\times\mathcal Y$, set $h_\xi=\mathcal R_0\xi$. Lax--Milgram on $\ker D_0$ gives a unique $k\in\ker D_0$ satisfying
\[
 \ip{A_0k}{z}_{\mathcal H}
 =\ip{f-A_0h_\xi}{z}_{\mathcal H}
 \qquad(z\in\ker D_0).
\]
The coercivity constant in \eqref{eq3.27} yields
\begin{equation}\label{eq3.31}
 \norm{k}_{\mathcal H}
 \le c_0^{-1}\left(
  \norm{f}_{\mathcal H}
  +\norm{A_0}\norm{\mathcal R_0}\norm{\xi}_{\mathcal Y}
  \right).
\end{equation}
For $h=h_\xi+k$, the vector $A_0h-f$ annihilates $\ker D_0$. Since $\operatorname{ran}D_0^*=(\ker D_0)^\perp$, there is a unique $\mu\in\mathcal Y^*$ with $A_0h-f=D_0^*\mu$. The map $D_0^*$ is injective on the finite-dimensional space $\mathcal Y^*$, so $\norm{\mu}\le C_D\norm{D_0^*\mu}_{\mathcal H}$. Together with \eqref{eq3.31}, this gives
\[
 \norm{h}_{\mathcal H}+\norm{\mu}_{\mathcal Y^*}
 \le C\left(\norm{f}_{\mathcal H}+\norm{\xi}_{\mathcal Y}\right).
\]
Thus \eqref{eq3.30} is an isomorphism.

The analytic implicit function theorem produces a unique analytic branch of constrained critical points $(u(\tau),\lambda(\tau))$ for $|\tau|$ small, with
\begin{equation}\label{eq3.32}
 u(\tau)=u_0+\tau u_2+O_{\mathcal H}(\tau^2).
\end{equation}
For $\tau=r^2\ge0$, the control $ru(r^2)$ on $G$ ends at $\exp(rX+r^2Z)$ and has energy $r^2E(u(r^2))=O(r^2)$.

We next identify the global minimizers with this critical branch. Properness of $(G,d_{\sr})$ gives a minimizing control $a_r$ on $G$ for the same endpoint. Set $v_r=a_r/r$. Since the local branch is an admissible competitor,
\begin{equation}\label{eq3.33}
 E(v_r)\le E(u(r^2)),
\end{equation}
so $(v_r)$ is bounded in $\mathcal H$. In particular, $\|rv_r\|_{L^1}=O(r)$, and \Cref{Lem3.1} places its endpoint in the logarithmic chart used in \eqref{eq3.24}. The exact endpoint identity is therefore equivalent to
\begin{equation}\label{eq3.34}
 F(r^2,v_r)=y.
\end{equation}

We claim that every choice of global minimizers satisfies
\begin{equation}\label{eq3.35}
 v_r\longrightarrow u_0\qquad\text{strongly in }\mathcal H.
\end{equation}
Let $r_n\downarrow0$. After passing to a subsequence, $v_{r_n}\rightharpoonup v$ in $\mathcal H$. Put $x_n(t)=\int_0^t v_{r_n}(s)\,\mathrm ds$ and $x(t)=\int_0^t v(s)\,\mathrm ds$. The Volterra operator is compact from $L^2$ to $C^0$, hence
\begin{equation}\label{eq3.36}
 x_n\longrightarrow x\qquad\text{in }C^0([0,1];\pLie).
\end{equation}
The horizontal tangent endpoint is weakly continuous. For the vertical endpoint,
\begin{align*}
 Z(v_{r_n})-Z(v)
 &=\frac12\int_0^1B((x_n-x)\wedge v_{r_n})\,\mathrm dt\\
 &\quad+\frac12\int_0^1B(x\wedge(v_{r_n}-v))\,\mathrm dt.
\end{align*}
The first integral tends to zero by \eqref{eq3.36} and the uniform $L^2$ bound. The second tends to zero by weak convergence. Thus $F_0$ is weakly sequentially continuous on bounded sets. From \eqref{eq3.25} and \eqref{eq3.34}, $F_0(v_{r_n})\to y$, and therefore $F_0(v)=y$.

Weak lower semicontinuity, uniqueness of the tangent minimizer, and \eqref{eq3.33} give
\[
 E(u_0)\le E(v)\le\liminf_{n\to\infty}E(v_{r_n})
 \le\limsup_{n\to\infty}E(v_{r_n})
 \le\lim_{n\to\infty}E(u(r_n^2))=E(u_0).
\]
Hence $v=u_0$ and the norms converge, so the convergence is strong. A contradiction argument shows that the same conclusion holds for every family of global minimizers, not merely for a selected subsequence. This proves \eqref{eq3.35}.

We identify these minimizers with the local multiplier branch. The $C^1$ estimate in \eqref{eq3.25} and \eqref{eq3.35} imply
\[
 D_uF(r^2,v_r)\longrightarrow D_0
 \qquad\text{in operator norm}.
\]
Set $D_r=D_uF(r^2,v_r)$. For small $r$, the operator
\[
 \mathcal R_r
 =\mathcal R_0\bigl(I+(D_r-D_0)\mathcal R_0\bigr)^{-1}
\]
is a right inverse of $D_r$, and $\sup_r\norm{\mathcal R_r}<\infty$. The Hilbert space Lagrange multiplier theorem gives a unique $\mu_r\in\mathcal Y^*$ such that
\[
 DE(v_r)=D_r^*\mu_r.
\]
For $\xi\in\mathcal Y$,
\[
 |\mu_r(\xi)|
 =|DE(v_r)[\mathcal R_r\xi]|
 \le2\norm{v_r}_{\mathcal H}\norm{\mathcal R_r}\norm{\xi}_{\mathcal Y},
\]
so $(\mu_r)$ is bounded. Every limit point $\mu$ satisfies $DE(u_0)=D_0^*\mu$, and uniqueness in \eqref{eq3.26} gives $\mu=\lambda_0$. Hence $(v_r,\mu_r)\to(u_0,\lambda_0)$. The local uniqueness supplied by the implicit function theorem now gives
\[
 v_r=u(r^2),\qquad \mu_r=\lambda(r^2)
\]
for all sufficiently small $r$. Thus the minimizing control on $G$ is unique.

Finally, the coefficient of $\tau$ in $F(\tau,u(\tau))=y$ is
\[
 D_0u_2+F_2(u_0)=0.
\]
Using \eqref{eq3.26} and \eqref{eq3.32},
\begin{align*}
 E(u(\tau))
 &=E(u_0)+2\tau\ip{u_0}{u_2}_{\mathcal H}+O(\tau^2)\\
 &=E(u_0)+\tau\lambda_0(D_0u_2)+O(\tau^2)\\
 &=E(u_0)-\tau\lambda_0(F_2(u_0))+O(\tau^2).
\end{align*}
Multiplication by $r^2$, with $\tau=r^2$, proves \eqref{eq3.28} and \eqref{eq3.29}.
\end{proof}

\begin{proposition}\label{Prop3.8}
Let $X\in\pLie$ and assume that the map
\[
 \pLie\longrightarrow\kLie,\qquad Y\longmapsto[X,Y],
\]
is onto. Then the conditions in \Cref{Ass3.6} hold at $(X,0)$. Its unique tangent minimizing control is $u_0(t)=X$, and the coefficient in \eqref{eq3.29} is $c_4(X,0)=0$.
\end{proposition}

\begin{proof}
For every control with horizontal endpoint $X$, the Cauchy--Schwarz inequality gives
\[
 E(u)\ge\left|\int_0^1u(t)\,\mathrm dt\right|^2=|X|^2.
\]
Equality holds only for the constant control $u_0(t)=X$, whose vertical endpoint is zero. Thus $u_0$ is the unique minimizer. A direct differentiation of the tangent endpoint map gives
\[
 DF_0(u_0)h=
 \left(\int_0^1h(t)\,\mathrm dt,
 \frac12\int_0^1(2t-1)[X,h(t)]\,\mathrm dt\right).
\]
Given $(P,Z)\in\pLie\oplus\kLie$, choose $Y\in\pLie$ with $[X,Y]=Z$ and set $h(t)=P+6(2t-1)Y$. Then $DF_0(u_0)h=(P,Z)$, so the differential is onto. The multiplier is $\lambda_0(P,Z)=2\ip{X}{P}$. Its vertical component is zero, and hence the constrained Hessian in \eqref{eq3.27} is $H_0(h,k)=2\ip{h}{k}_{\mathcal H}$. Finally, all brackets entering $Q(u_0)$ and $R(u_0)$ vanish, so $c_4(X,0)=0$.
\end{proof}

Surjectivity of $Y\mapsto[X,Y]$ is an open condition on $X$. For the standard symmetric pair in $\SU(2)$, this map is onto for every $X\ne0$. Thus \Cref{Prop3.8} applies to every nonzero horizontal datum $(X,0)$ in that example.

The estimate in \Cref{Thm3.4} holds uniformly on compact sets of tangent endpoints, including points where one or more of the conditions in \Cref{Ass3.6} fail. Under all three conditions, \Cref{Thm3.7} determines the fourth-order coefficient and improves the remainder to $O(r^6)$.

\subsection{Sharpness at vertical cut points}

The following example on $\SU(2)$ shows that the $O(r^4)$ remainder in \Cref{Thm3.4} is sharp.

\begin{proposition}
\label{Prop3.9}
Equip $\mathfrak{su}(2)$ with $\ip{A}{B}=-\Tr(AB)$, let
\[
 k=\diag(i/2,-i/2),\qquad
 \kLie=\R k,\qquad \pLie=\kLie^\perp,
\]
and let $d_{\sr}$ be generated by $\pLie$. For fixed $z\in\R$ and all sufficiently small $r>0$,
\begin{equation}\label{eq3.37}
 d_{\sr}\!\left(e,\exp(r^2zk)\right)^2
 =2\pi|z|r^2-\frac12z^2r^4.
\end{equation}
For every fixed $z\ne0$, the coefficient of $r^4$ in \eqref{eq3.37} is nonzero. Hence the remainder in \Cref{Thm3.4} cannot be replaced by $o(r^4)$ uniformly on a compact set of tangent data containing such vertical points.
\end{proposition}

\begin{proof}
The exact vertical formula is given in \cite[Theorem~2]{BoscainRossi2008} and \cite[Theorem~2]{BerestovskiiZubareva2016}. In the orthonormal normalization used there, it is
\[
 d_{\mathrm{BZ}}(e,\diag(e^{i\theta/2},e^{-i\theta/2}))^2
 =4\frac{|\theta|}{2}\left(2\pi-\frac{|\theta|}{2}\right)
\]
for $\theta$ near zero. Their horizontal inner product is $2(-\Tr(AB))$, so squared distances in the present normalization are one half of the displayed value. Taking $\theta=r^2z$ proves \eqref{eq3.37}. Its leading coefficient also gives $d_0(e,(0,zk))^2=2\pi|z|$, as required by \Cref{Thm3.4}.
\end{proof}

\section{Vertical squared distance and block matrix formulas}
\label{Sec4}

The comparison theorems of \Cref{Sec3} reduce the problem near the identity to the squared tangent distance. We first recall the formula for the vertical squared distance proved in \cite{BauerHanYang2026} and then evaluate it for block subgroups of $\SU(N)$. The free group identity also agrees with the Fourier isoperimetric formula of Rizzi and Serres \cite[Proposition~15]{RizziSerres2017}.

Let $V$ be a real Euclidean space of dimension $n$. For $\Omega\in\Lambda^2V$, define the skew-symmetric operator $J_\Omega$ by
\[
 \ip{J_\Omega v}{w}=\ip{\Omega}{v\wedge w},
 \qquad v,w\in V,
\]
and write
\[
 \sigma_1(\Omega)\ge\cdots\ge
 \sigma_{\lfloor n/2\rfloor}(\Omega)\ge0
\]
for its nonnegative normal form coefficients. Set
\begin{equation}\label{eq4.1}
 \Psi_V(\Omega)
 =4\pi\sum_{j=1}^{\lfloor n/2\rfloor}j\,\sigma_j(\Omega).
\end{equation}
For the Lie pair of \Cref{Sec2}, recall the surjective projected bracket
\[
 B:\Lambda^2\pLie\longrightarrow\lLie,
 \qquad B(X\wedge Y)=\pr_{\lLie}[X,Y],
\]
and put $\mathfrak F_{G,L}(Z)=d_0(e,(0,Z))^2$.

\begin{proposition}\label{Prop4.1}
For every $Z\in\lLie$, with $\Psi_{\pLie}$ defined by \eqref{eq4.1},
\begin{equation}\label{eq4.2}
 \mathfrak F_{G,L}(Z):=d_0(e,(0,Z))^2
 =\min_{B\Omega=Z}\Psi_{\pLie}(\Omega)
 =\min_{B\Omega=Z}
 4\pi\sum_{j=1}^{\lfloor\dim\pLie/2\rfloor}
 j\,\sigma_j(\Omega).
\end{equation}
The minimum is attained, and every minimizing lift gives a vertical minimizing loop using at most $\lfloor\dim\pLie/2\rfloor$ positive Fourier frequencies. Moreover, $\mathfrak F_{G,L}$ is positively homogeneous of degree one, globally Lipschitz, and equivalent to the Euclidean norm on $\lLie$: there are $c,C,L>0$, depending only on the pair, such that
\begin{equation}\label{eq4.3}
 c|Z|\le\mathfrak F_{G,L}(Z)\le C|Z|,
 \qquad
 |\mathfrak F_{G,L}(Z)-\mathfrak F_{G,L}(W)|\le L|Z-W|.
\end{equation}
\end{proposition}

\begin{proof}
The free two-step identity
\[
 d_{\mathrm{free}}(e,(0,\Omega))^2=\Psi_{\pLie}(\Omega)
\]
is \cite[Theorem~4.4]{BauerHanYang2026}. For the projected bracket in \eqref{eq2.3}, the map
\[
 q_B:\mathbb F_2(\pLie)\longrightarrow\mathbb N_{G,L},
 \qquad q_B(X,\Omega)=(X,B\Omega),
\]
is a central Carnot group quotient and a submetry. See \cite[Proposition~4.6, Corollaries~4.8--4.9, and Remark~4.13]{BauerHanYang2026}. Taking the infimum over the fiber $B^{-1}(Z)$ proves \eqref{eq4.2}. The statement about minimizing loops follows from the corresponding statement on the free group. The fiber is nonempty and closed. Since $\Psi_{\pLie}(\Omega)\ge4\pi|\Omega|$, the minimum is attained. Finally, the proof of \cite[Proposition~4.1]{BauerHanYang2026} uses only surjectivity of $B$. The Hoffman--Wielandt inequality gives the Lipschitz bound for $\Psi_{\pLie}$, and comparison along a linear right inverse of $B$ gives \eqref{eq4.3}. Homogeneity follows directly from \eqref{eq4.2}.
\end{proof}

For block subgroups, the affine fibers in \eqref{eq4.2} are described by linear constraints on the partial traces of Hermitian matrices.

\subsection{Block subgroups}\label{Subsec4.1}
The same block family appears in \cite[Section~6]{BauerHanYang2026}, where the horizontal complement is computed and the two-step condition is proved. We retain only the notation needed below.

Equip $\mathfrak{su}(N)$ with $\ip{X}{Y}=-\Tr(XY)$, fix a partition
\[
 N=n_1+\cdots+n_s,\qquad s\ge2,\qquad n_a\ge1,
\]
and put $\mathbf n=(n_1,\ldots,n_s)$. Define the traceless block scalar algebra and its connected torus by
\begin{align*}
 \mathfrak z_{\mathbf n}
 &=\left\{\diag(i\beta_1I_{n_1},\ldots,i\beta_sI_{n_s}):
 \beta_a\in\R,\ \sum_{a=1}^sn_a\beta_a=0\right\},\\
 Z_{\mathbf n}&:=\exp(\mathfrak z_{\mathbf n}).
\end{align*}
Let $T_0\subset Z_{\mathbf n}$ be a closed connected subtorus with Lie algebra $\mathfrak t_0$, and set
\[
 L_{\mathbf n,T_0}
 =\left(\prod_{a=1}^s\SU(n_a)\right)T_0\subset\SU(N).
\]
Put
\[
 V_{\mathbf n}
 =\left\{\beta\in\R^s:\sum_{a=1}^sn_a\beta_a=0\right\},
 \qquad
 \ip{\beta}{\gamma}_{\mathbf n}
 =\sum_{a=1}^sn_a\beta_a\gamma_a,
\]
and let $V_0\subset V_{\mathbf n}$ correspond to $\mathfrak t_0$. This family includes $(\SU(N),\SU(N-1))$, obtained from $\mathbf n=(N-1,1)$ and $T_0=\{e\}$, and the full block Levi subgroup
\[
 S(\U(n_1)\times\cdots\times\U(n_s)),
\]
obtained by taking $T_0=Z_{\mathbf n}$.

Let $\lLie=\Lie(L_{\mathbf n,T_0})$ and $\pLie=\lLie^\perp$. As shown in \cite[Section~6]{BauerHanYang2026}, every element of $\pLie$ has the unique form
\begin{equation}\label{eq4.4}
 X(\beta,Z)=
 \begin{pmatrix}
  i\beta_1I_{n_1}&Z_{12}&\cdots&Z_{1s}\\
  -Z_{12}^*&i\beta_2I_{n_2}&\cdots&Z_{2s}\\
  \vdots&\vdots&\ddots&\vdots\\
  -Z_{1s}^*&-Z_{2s}^*&\cdots&i\beta_sI_{n_s}
 \end{pmatrix},
 \qquad \beta\in V_0^\perp,
\end{equation}
where $Z_{ab}\in M_{n_a\times n_b}(\C)$, and
\begin{equation}\label{eq4.5}
 |X(\beta,Z)|^2
 =\sum_{a=1}^sn_a\beta_a^2
  +2\sum_{a<b}\norm{Z_{ab}}_{\mathrm{HS}}^2.
\end{equation}
Moreover,
\[
 \mathfrak{su}(N)=\pLie+[\pLie,\pLie]
\]
by \cite[Proposition~6.1]{BauerHanYang2026}. Hence the quotient formula \eqref{eq4.2} applies to every block pair above.

We use Hermitian coordinates for the vertical layer: a Hermitian matrix $H$ represents the Lie algebra element $iH$. Put
\begin{align*}
 \mathcal W_{\mathbf n}
 &=\left\{\diag(U_1,\ldots,U_s):U_a=U_a^*,\
  \sum_a\Tr U_a=0\right\},
 \\
 \mathcal V_{\mathbf n,T_0}^{\mathrm H}
 &=\left\{H\in\mathcal W_{\mathbf n}:
  \left(\frac{\Tr H_1}{n_1},\ldots,\frac{\Tr H_s}{n_s}\right)\in V_0
  \right\}.
\end{align*}
Write $\mathbf D(\beta)=\diag(\beta_1I_{n_1},\ldots,\beta_sI_{n_s})$. There is an orthogonal decomposition
\begin{equation*}
 \mathcal W_{\mathbf n}
 =\mathbf D(V_0^\perp)\oplus
  \mathcal V_{\mathbf n,T_0}^{\mathrm H}.
\end{equation*}
Let $\pi_{\mathbf n,T_0}$ denote the projection onto the second summand.

For a based off-diagonal block path $Z=(Z_{ab})_{a<b}$, define $\mathcal U(Z)=\diag(\mathcal U_1(Z),\ldots,\mathcal U_s(Z))$ by
\begin{align}
 \mathcal U_a(Z)=-\frac{i}{2}\Biggl(
 &\sum_{b>a}\int_0^1
  (\dot Z_{ab}Z_{ab}^*-Z_{ab}\dot Z_{ab}^*)\,\mathrm dt
 \notag\\[-0.2em]
 &+\sum_{c<a}\int_0^1
  (\dot Z_{ca}^*Z_{ca}-Z_{ca}^*\dot Z_{ca})\,\mathrm dt
 \Biggr),
 \label{eq4.6}
\end{align}
and
\begin{equation}\label{eq4.7}
 \mathcal E(Z)=2\sum_{a<b}\int_0^1
  \norm{\dot Z_{ab}(t)}_{\mathrm{HS}}^2\,\mathrm dt.
\end{equation}
Set
\begin{equation}\label{eq4.8}
 \mathfrak I_{\mathbf n,T_0}(H)
 =\inf\left\{\mathcal E(Z):
 Z_{ab}(0)=Z_{ab}(1)=0,\
 \pi_{\mathbf n,T_0}\mathcal U(Z)=H\right\}.
\end{equation}

\begin{lemma}
\label{Lem4.2}
For $H\in\mathcal V_{\mathbf n,T_0}^{\mathrm H}$,
\begin{equation}\label{eq4.9}
 \mathfrak F_{G,L}(iH)=\mathfrak I_{\mathbf n,T_0}(H).
\end{equation}
\end{lemma}

\begin{proof}
A tangent horizontal loop may have a block scalar component together with the off-diagonal components in \eqref{eq4.4}. Brackets involving a block scalar horizontal component are off-diagonal and hence have zero $\lLie$-projection. Two scalar components commute. Brackets from two different block edges are again off-diagonal. Thus only an edge paired with itself contributes to the vertical endpoint.

For a path supported on one edge $a<b$, direct multiplication and integration give
\begin{align*}
 \frac12\int_0^1
 [X(0,Z(t)),X(0,\dot Z(t))]_{aa}\,\mathrm dt
 &=\frac12\int_0^1(\dot ZZ^*-Z\dot Z^*)\,\mathrm dt\\
 &=i\,\mathcal U_a(Z),
\end{align*}
and the lower diagonal block gives the incoming term in \eqref{eq4.6}. Therefore the projected tangent area is precisely $i\pi_{\mathbf n,T_0}\mathcal U(Z)$. The energy is \eqref{eq4.7} by \eqref{eq4.5}. Setting the scalar component equal to zero preserves the vertical constraint and cannot increase energy, proving \eqref{eq4.9}.
\end{proof}

\subsection{A finite-dimensional matrix formula}\label{Subsec4.2}

For $a<b$, let
\[
 \mathscr H_{ab}=M_{n_a\times n_b}(\C),
 \qquad d_{ab}=n_an_b,
\]
with its Hilbert--Schmidt inner product. For $A\in\mathscr H_{ab}$, the notation $|A\rangle\langle A|$ denotes the rank-one operator
\[
 (|A\rangle\langle A|)Z=\Tr(A^*Z)\,A,
 \qquad Z\in\mathscr H_{ab}.
\]
If $S\in\Herm(\mathscr H_{ab})$ has entries $S_{i\mu,j\nu}$ in the standard matrix unit basis, define the linear maps $\mathsf L_{ab}$ and $\mathsf R_{ab}$ by
\begin{align}
 (\mathsf L_{ab}S)_{ij}
 &=\sum_{\mu=1}^{n_b}S_{i\mu,j\mu},\notag\\
 (\mathsf R_{ab}S)_{\mu\nu}
 &=\sum_{i=1}^{n_a}S_{i\nu,i\mu}.
 \label{eq4.10}
\end{align}
In the chosen matrix unit coordinates, $\mathsf L_{ab}$ is a partial trace and $\mathsf R_{ab}$ is the transpose of the other partial trace. The reversed indices in the second formula ensure that
\begin{equation}\label{eq4.11}
 \mathsf L_{ab}(|A\rangle\langle A|)=AA^*,
 \qquad
 \mathsf R_{ab}(|A\rangle\langle A|)=A^*A.
\end{equation}
Both resulting matrices have trace $\Tr S$.

On
\[
 \mathscr S_{\mathbf n}
 =\bigoplus_{a<b}\Herm(\mathscr H_{ab}),
 \qquad
 \norm{S}_{\mathscr S_{\mathbf n}}^2
 =\sum_{a<b}\norm{S_{ab}}_{\mathrm{HS}}^2,
\]
define the linear map $\mathscr A_{\mathbf n}(S) =\diag(\mathscr A_1(S),\ldots,\mathscr A_s(S))$ by
\begin{equation*}
 \mathscr A_a(S)
 =\sum_{b>a}\mathsf L_{ab}(S_{ab})
  -\sum_{c<a}\mathsf R_{ca}(S_{ca}),
\end{equation*}
and set
\begin{equation*}
 \mathscr Q_{\mathbf n,T_0}
 =\pi_{\mathbf n,T_0}\circ\mathscr A_{\mathbf n}.
\end{equation*}
The equality of these traces gives $\mathscr A_{\mathbf n}(S)\in\mathcal W_{\mathbf n}$.

\begin{lemma}
\label{Lem4.3}
The map $\mathscr A_{\mathbf n}$ is onto $\mathcal W_{\mathbf n}$. Hence $\mathscr Q_{\mathbf n,T_0}$ is onto $\mathcal V_{\mathbf n,T_0}^{\mathrm H}$.
\end{lemma}

\begin{proof}
For unit vectors $u\in\C^{n_a}$, $v\in\C^{n_b}$, set $A=uv^*$. By \eqref{eq4.11}, the image of $|A\rangle\langle A|$ has $uu^*$ in block $a$, $-vv^*$ in block $b$, and zero elsewhere. Differences with $v$ fixed span the traceless Hermitian matrices in block $a$. Summing over orthonormal bases of $u$ and $v$ gives $n_bI_{n_a}$ in block $a$ and $-n_aI_{n_b}$ in block $b$. These edge differences span the traceless block scalar space. The two families together span $\mathcal W_{\mathbf n}$.
\end{proof}

For $S\in\Herm(d)$, write $\lambda_1(S)\ge\cdots\ge\lambda_d(S)$ for its eigenvalues. As in \cite{BauerHanYang2026}, let
\[
 \lambda_1^+(S)\ge\cdots\ge\lambda_{p(S)}^+(S)>0,
 \qquad
 \lambda_1^-(S)\ge\cdots\ge\lambda_{q(S)}^-(S)>0
\]
be the positive eigenvalues and the absolute values of the negative eigenvalues, respectively. Define
\begin{equation}\label{eq4.12}
 \Phi_d(S)
 :=4\pi\left(
 \sum_{j=1}^{p(S)}j\lambda_j^+(S)
 +\sum_{j=1}^{q(S)}j\lambda_j^-(S)
 \right).
\end{equation}
An empty sum is zero. Equivalently,
\[
 \Phi_d(S)=4\pi\sum_{j=1}^d j
 \left((\lambda_j(S))_+
 +(-\lambda_{d+1-j}(S))_+\right).
\]

\begin{lemma}
\label{Lem4.4}
For $S\in\Herm(\mathscr H)$, $\dim_{\C}\mathscr H=d$,
\begin{equation}\label{eq4.13}
 \frac{\Phi_d(S)}{4\pi}
 =\min\Biggl\{
 \sum_{k\ge1}k\bigl(\norm{A_k^+}^2+\norm{A_k^-}^2\bigr):
 \notag
 S=\sum_{k\ge1}|A_k^+\rangle\langle A_k^+|
  -\sum_{k\ge1}|A_k^-\rangle\langle A_k^-|
 \Biggr\}.
\end{equation}
The series are understood in trace norm, and a minimizer uses no index larger than $d$.
\end{lemma}

\begin{proof}
First let $M\ge0$, with eigenvalues $\mu_1\ge\cdots\ge\mu_d\ge0$. We claim
\begin{equation}\label{eq4.14}
 \min\left\{\sum_{k\ge1}k\Tr X_k:
 M=\sum_{k\ge1}X_k,\ X_k\ge0,\ \rank X_k\le1\right\}
 =\sum_{j=1}^d j\mu_j,
\end{equation}
where the positive series converges in trace norm. The spectral decomposition gives the upper bound. For any admissible decomposition, $P_m=\sum_{k\le m}X_k$ satisfies $0\le P_m\le M$ and $\rank P_m\le m$. Ky Fan's principle gives
\[
 \Tr P_m\le\sum_{j\le m}\mu_j,
 \qquad
 \sum_{k>m}\Tr X_k\ge\sum_{j>m}\mu_j.
\]
Summing over $m\ge0$, which is legitimate for the nonnegative trace series, proves \eqref{eq4.14}.

For a signed decomposition, define the positive operators
\[
 P=\sum_k|A_k^+\rangle\langle A_k^+|,
 \qquad
 Q=\sum_k|A_k^-\rangle\langle A_k^-|.
\]
Both series converge in trace norm. Then $S=P-Q$, so $P\ge S$ and $Q\ge-S$. Let $p$ be the number of positive eigenvalues of $S$. Weyl monotonicity gives
\begin{equation}\label{eq4.15}
 \lambda_j(P)\ge\lambda_j(S)>0,
 \qquad 1\le j\le p.
\end{equation}
If $q$ is the number of negative eigenvalues of $S$, then
\begin{equation}\label{eq4.16}
 \lambda_j(Q)\ge\lambda_j(-S)=-\lambda_{d+1-j}(S)>0,
 \qquad 1\le j\le q.
\end{equation}
Applying \eqref{eq4.14} to $P$ and $Q$, and then using \eqref{eq4.15}--\eqref{eq4.16}, gives the lower bound in \eqref{eq4.13}. Equality is obtained by placing the positive spectral vectors of $S$, in decreasing eigenvalue order, into the positive modes $1,\ldots,p$, and the negative spectral vectors, in decreasing absolute eigenvalue order, into the negative modes $1,\ldots,q$. Thus a minimizer uses no index larger than $d$.
\end{proof}

\begin{lemma}
\label{Lem4.5}
Let $Z\in H^1([0,1];M_{m\times n}(\C))$ be based and write
\[
 Z(t)=\sum_{k\in\Z}Z_ke^{2\pi i kt}.
\]
Then
\begin{align}
 U_L(Z)
 &:=-\frac{i}{2}\int_0^1(\dot ZZ^*-Z\dot Z^*)\,\mathrm dt
  =2\pi\sum_{k\in\Z}kZ_kZ_k^*,
 \label{eq4.17}\\
 U_R(Z)
 &:=-\frac{i}{2}\int_0^1(\dot Z^*Z-Z^*\dot Z)\,\mathrm dt
  =-2\pi\sum_{k\in\Z}kZ_k^*Z_k,\notag\\
 E(Z)
 &:=2\int_0^1\norm{\dot Z}_{\mathrm{HS}}^2\,\mathrm dt
  =8\pi^2\sum_{k\in\Z}k^2\norm{Z_k}_{\mathrm{HS}}^2.
 \label{eq4.18}
\end{align}
For $k\ge1$, put
\[
 A_k^+=\sqrt{2\pi k}\,Z_k,\qquad
 A_k^-=\sqrt{2\pi k}\,Z_{-k},
\]
and
\begin{equation}\label{eq4.19}
 S_Z=\sum_{k\ge1}|A_k^+\rangle\langle A_k^+|
  -\sum_{k\ge1}|A_k^-\rangle\langle A_k^-|.
\end{equation}
Then the series converges in trace norm,
\begin{equation*}
 U_L(Z)=\mathsf L(S_Z),\qquad
 U_R(Z)=-\mathsf R(S_Z),
\end{equation*}
and
\begin{equation}\label{eq4.20}
 E(Z)\ge\Phi_{mn}(S_Z).
\end{equation}
Conversely, every $S\in\Herm(M_{m\times n}(\C))$ is realized by a based loop using only frequencies $|k|\le mn$ and having energy $\Phi_{mn}(S)$.
\end{lemma}

\begin{proof}
The endpoint values agree, so the based path defines an $H^1$-loop on the circle. Its Fourier series converges in $H^1$ and uniformly. Applying Fourier orthogonality to partial sums and passing to the limit gives \eqref{eq4.17}--\eqref{eq4.18}. Moreover,
\[
 \sum_{k\ge1}\left(\norm{A_k^+}_{\mathrm{HS}}^2
  +\norm{A_k^-}_{\mathrm{HS}}^2\right)
 =2\pi\sum_{k\ne0}|k|\norm{Z_k}_{\mathrm{HS}}^2<\infty.
\]
Hence both positive series in \eqref{eq4.19} converge in trace norm. The identities for $\mathsf L(S_Z)$ and $\mathsf R(S_Z)$ follow from \eqref{eq4.11}, and
\[
 E(Z)=4\pi\sum_{k\ge1}k
  \left(\norm{A_k^+}_{\mathrm{HS}}^2
  +\norm{A_k^-}_{\mathrm{HS}}^2\right).
\]
\Cref{Lem4.4} proves \eqref{eq4.20}. Notice that, for a fixed $k>0$, the coefficients $Z_k$ and $Z_{-k}$ are independent. Thus the same frequency magnitude can carry one positive and one negative rank-one term.

For the converse, let
\[
 S=\sum_{j=1}^p\lambda_j^+
  |U_j^+\rangle\langle U_j^+|
  -\sum_{j=1}^q\lambda_j^-
  |U_j^-\rangle\langle U_j^-|
\]
be an orthonormal spectral decomposition, with both eigenvalue lists in decreasing order. The loop
\begin{equation}\label{eq4.21}
 Z_S(t)=
 \sum_{j=1}^p\sqrt{\frac{\lambda_j^+}{2\pi j}}\,
  (e^{2\pi i jt}-1)U_j^+
 +\sum_{j=1}^q\sqrt{\frac{\lambda_j^-}{2\pi j}}\,
  (e^{-2\pi i jt}-1)U_j^-
\end{equation}
is based, realizes the matrix $S$, and has energy $\Phi_{mn}(S)$.
\end{proof}

For $H\in\mathcal V_{\mathbf n,T_0}^{\mathrm H}$, set
\begin{equation}\label{eq4.22}
 \mathfrak F_{\mathbf n,T_0}(H)
 =\min\left\{\sum_{a<b}\Phi_{d_{ab}}(S_{ab}):
 S\in\mathscr S_{\mathbf n},\
 \mathscr Q_{\mathbf n,T_0}(S)=H\right\}.
\end{equation}

\begin{theorem}
\label{Thm4.6}
For every $H\in\mathcal V_{\mathbf n,T_0}^{\mathrm H}$,
\begin{equation}\label{eq4.23}
 d_0(e,(0,iH))^2=\mathfrak F_{\mathbf n,T_0}(H).
\end{equation}
The minimum in \eqref{eq4.22} is attained. If $S=(S_{ab})$ is a minimizer, applying \eqref{eq4.21} to each pair $a<b$ gives a tangent minimizing loop involving no frequency larger than $n_an_b$ in the corresponding component.

The coefficient $\mathfrak F_{\mathbf n,T_0}$ is globally Lipschitz and positively homogeneous of degree one. There are constants $c,C>0$ such that
\begin{equation}\label{eq4.24}
 c\norm{H}_{\mathrm{HS}}\le\mathfrak F_{\mathbf n,T_0}(H)
 \le C\norm{H}_{\mathrm{HS}},
\end{equation}
where $\norm{\cdot}_{\mathrm{HS}}$ is the Hilbert--Schmidt norm on the finite-dimensional Hermitian block diagonal space $\mathcal V_{\mathbf n,T_0}^{\mathrm H}$.
\end{theorem}

\begin{proof}
We first identify the variational problems and prove attainment. Let $Z=(Z_{ab})$ be feasible in \eqref{eq4.8}, and form $S_{ab}=S_{Z_{ab}}$ using \eqref{eq4.19}. The identities for $\mathsf L_{ab}$ and $\mathsf R_{ab}$ show, block by block, that
\[
 \mathscr A_{\mathbf n}(S)=\mathcal U(Z),\qquad
 \mathscr Q_{\mathbf n,T_0}(S)=H.
\]
By \eqref{eq4.20},
\[
 \sum_{a<b}\Phi_{d_{ab}}(S_{ab})\le\mathcal E(Z).
\]
This proves $\mathfrak F_{\mathbf n,T_0}\le\mathfrak I_{\mathbf n,T_0}$. Conversely, any feasible family $S$ produces, by \eqref{eq4.21}, a based block loop with $\mathcal U(Z)=\mathscr A_{\mathbf n}(S)$ and energy equal to the value in \eqref{eq4.22}. Thus the reverse inequality holds. Together with \Cref{Lem4.2}, this gives
\[
 \mathfrak F_{\mathbf n,T_0}(H)
 =\mathfrak I_{\mathbf n,T_0}(H)
 =\mathfrak F_{G,L}(iH)
 =d_0(e,(0,iH))^2,
\]
which proves \eqref{eq4.23} and the statement about the minimizing loop.

Feasibility follows from \Cref{Lem4.3}. Since all weights in $\Phi_d$ are at least one,
\begin{equation*}
 \Phi_d(S)\ge4\pi\Tr|S|.
\end{equation*}
The constraint in \eqref{eq4.22} is a closed affine subspace. Every minimizing sequence is therefore bounded and has a convergent subsequence. This proves attainment.

For regularity and norm equivalence, rewrite \eqref{eq4.12} as
\[
 \Phi_d(S)=4\pi\sum_{j=1}^d
 \left[j(\lambda_j(S))_+
 +(d+1-j)(-\lambda_j(S))_+\right].
\]
The scalar summands have slopes bounded by $d$. The Hoffman--Wielandt inequality therefore gives
\begin{equation}\label{eq4.25}
 |\Phi_d(S)-\Phi_d(T)|
 \le4\pi d\sqrt d\,\norm{S-T}_{\mathrm{HS}}.
\end{equation}
Consequently, the total objective $J(S)=\sum_{a<b}\Phi_{d_{ab}}(S_{ab})$ is Lipschitz with respect to $\norm{\cdot}_{\mathscr S_{\mathbf n}}$. The required eigenvalue inequalities may be found in \cite[Chapters~III,~VI]{Bhatia1997}. Choose a linear right inverse $\mathscr R$ of $\mathscr Q_{\mathbf n,T_0}$ and fix $C_R$ such that
\begin{equation}\label{eq4.26}
 \norm{\mathscr RH}_{\mathscr S_{\mathbf n}}
 \le C_R\norm{H}_{\mathrm{HS}}.
\end{equation}
If $S_H$ minimizes at $H$, then $S_H+\mathscr R(H'-H)$ is feasible at $H'$. Applying the Lipschitz bound for $J$ and interchanging $H,H'$ proves the global Lipschitz estimate for $\mathfrak F_{\mathbf n,T_0}$. Positive homogeneity of degree one follows directly from the constraint and the homogeneity of $\Phi_d$.

Since $J(0)=0$, \eqref{eq4.25} and \eqref{eq4.26} give the upper inequality in \eqref{eq4.24}. For every feasible $S$,
\[
 \norm{H}_{\mathrm{HS}}
 \le\norm{\mathscr Q_{\mathbf n,T_0}}
  \norm{S}_{\mathscr S_{\mathbf n}},
 \qquad
 J(S)\ge4\pi\sum_{a<b}\norm{S_{ab}}_1
  \ge4\pi\norm{S}_{\mathscr S_{\mathbf n}}.
\]
Taking the minimum gives the lower inequality.
\end{proof}

\subsection{Explicit block families}\label{Subsec4.3}

The finite-dimensional problem can be evaluated explicitly for several block families. The first one is a scalar transportation problem.

\begin{proposition}
\label{Prop4.7}
Assume $\mathbf n=(1,\ldots,1)$, put
\[
 V=\left\{x\in\R^s:\sum_{a=1}^s x_a=0\right\},
\]
and let $P_0:V\to V_0$ be the orthogonal projection. Identify $H=\diag(h_1,\ldots,h_s)\in\mathcal V_{\mathbf n,T_0}^{\mathrm H}$ with $h\in V_0$. Orient the edge $a b$ from $a$ to $b$ when $a<b$, and define the incidence map
\begin{equation*}
 (\mathsf Bx)_a=\sum_{b>a}x_{ab}-\sum_{c<a}x_{ca}.
\end{equation*}
Then
\begin{align}
 \mathfrak F_{\mathbf n,T_0}(H)
 &=4\pi\min\left\{\sum_{a<b}|x_{ab}|:P_0\mathsf Bx=h\right\}
 \label{eq4.27}\\
 &=2\pi\min_{z\in V_0^\perp}\|h+z\|_{\ell^1}
 \label{eq4.28}\\
 &=4\pi\max\left\{
 \ip{h}{y}:y\in V_0,\ \max_a y_a-\min_a y_a\le1\right\}.
 \label{eq4.29}
\end{align}
If $T_0$ is the full block scalar torus, then $V_0=V$ and
\begin{equation}\label{eq4.30}
  \mathfrak F_{(1,\ldots,1),T_0}(\diag(h_1,\ldots,h_s))
 =2\pi\sum_{a=1}^s|h_a|.
\end{equation}
\end{proposition}

\begin{proof}
For a scalar edge, $S_{ab}=x_{ab}\in\R$, both partial traces equal $x_{ab}$, and $\Phi_1(x_{ab})=4\pi|x_{ab}|$. The incidence map is therefore $\mathsf Bx$, which proves \eqref{eq4.27}.

For every $g\in V$,
\begin{equation}\label{eq4.31}
 \min_{\mathsf Bx=g}\sum_{a<b}|x_{ab}|=\frac12\|g\|_{\ell^1}.
\end{equation}
Indeed, summing the components over $\{a:g_a>0\}$ gives the lower bound. A transportation plan from the positive masses $g_a$ to the negative masses $-g_b$, with signs adjusted to the fixed edge orientations, attains it. Since $P_0\mathsf Bx=h$ is equivalent to $\mathsf Bx=h+z$ for some $z\in V_0^\perp$, \eqref{eq4.31} gives \eqref{eq4.28}.

Linear programming duality in finite dimensions, applied to \eqref{eq4.27}, gives
\[
 \frac{\mathfrak F_{\mathbf n,T_0}(H)}{4\pi}
 =\max\left\{\ip{h}{y}:y\in V_0,
  \|\mathsf B^*y\|_{\ell^\infty}\le1\right\}.
\]
Here $(\mathsf B^*y)_{ab}=y_a-y_b$. Since the graph is complete, the dual constraint is the oscillation constraint in \eqref{eq4.29}. If $V_0=V$, then $V_0^\perp=\{0\}$, and \eqref{eq4.30} follows from \eqref{eq4.28}.
\end{proof}

The second family has only one edge. For the full block scalar torus, the one-column AIII formula
\[
 \mathfrak F_{(m,1),T_0}(\diag(A,-\Tr A))=\Phi_m(A)
\]
and its minimizing loop using finitely many Fourier frequencies were established in \cite[Theorem~6.2]{BauerHanYang2026}. We therefore focus on the different quotient obtained when the block scalar vertical direction is removed.

\begin{proposition}
\label{Prop4.8}
Let $\mathbf n=(m,1)$ and $T_0=\{e\}$. Every vertical matrix has the form
\[
 \widehat H_A=\diag(A,0),
 \qquad A\in\Herm(m),\qquad\Tr A=0.
\]
If $S\in\Herm(m)$ denotes the single edge matrix, then
\begin{equation}\label{eq4.32}
 \mathscr Q_{(m,1),\{e\}}(S)
 =\diag\left(S-\frac{\Tr S}{m}I_m,0\right).
\end{equation}
Writing $\lambda_1(A)\ge\cdots\ge\lambda_m(A)$ and $r_m=\lfloor m/2\rfloor$, one has
\begin{align}
 \mathfrak F_{(m,1),\{e\}}(\widehat H_A)
 &=\min_{t\in\R}\Phi_m(A+tI_m)
 \label{eq4.33}\\
 &=4\pi\sum_{j=1}^{r_m}j
  \bigl(\lambda_j(A)-\lambda_{m+1-j}(A)\bigr).
 \label{eq4.34}
\end{align}
The minimizing shifts are
\begin{equation}\label{eq4.35}
 \operatorname*{argmin}_{t\in\R}\Phi_m(A+tI_m)
 =\begin{cases}
 \{-\lambda_{r_m+1}(A)\},&m=2r_m+1,\\
 [-\lambda_{r_m}(A),-\lambda_{r_m+1}(A)],&m=2r_m.
 \end{cases}
\end{equation}
\end{proposition}

\begin{proof}
For a single edge, the image under $\mathscr A_{\mathbf n}$ is $\diag(S,-\Tr S)$. Put $\bar s=\Tr S/m$. The orthogonal splitting
\[
 \diag(S,-\Tr S)
 =\diag(S-\bar sI_m,0)+\diag(\bar sI_m,-m\bar s)
\]
shows that the projection onto the Hermitian vertical space is exactly \eqref{eq4.32}. Thus the feasible matrices are $S=A+tI_m$, which proves \eqref{eq4.33}.

Set $f(t)=\Phi_m(A+tI_m)/(4\pi)$. Reindexing its negative part gives
\[
 f(t)=\sum_{j=1}^m\left[
 j(\lambda_j+t)_+ +(m+1-j)(-\lambda_j-t)_+\right].
\]
On each open interval between consecutive breakpoints $-\lambda_j(A)$, if precisely the first $k$ shifted eigenvalues are positive, then
\[
 f'(t)=\frac{k(k+1)}2-\frac{(m-k)(m-k+1)}2
  =\frac{m+1}{2}(2k-m).
\]
The function is convex and continuous, and its subdifferential at a breakpoint is the interval between the adjacent slopes. For $m=2r_m+1$, the slope changes sign at $t=-\lambda_{r_m+1}$. For $m=2r_m$, it vanishes throughout $[-\lambda_{r_m},-\lambda_{r_m+1}]$. This proves \eqref{eq4.35}. Substitution of any minimizing shift, followed by pairing the upper and lower eigenvalues, yields \eqref{eq4.34}.
\end{proof}

The following result treats one matrix block together with several scalar blocks.

\begin{lemma}
\label{Lem4.9}
For $M\in\Herm(m)$ and $q\ge1$, define
\[
 \Psi_{m,q}(M)
 =\min\left\{\sum_{a=1}^q\Phi_m(S_a):
 S_a\in\Herm(m),\ \sum_{a=1}^qS_a=M\right\}.
\]
Then
\begin{equation}\label{eq4.36}
  \Psi_{m,q}(M)
 =4\pi\sum_{j=1}^m\left\lceil\frac jq\right\rceil
 \left((\lambda_j(M))_+
  +(-\lambda_{m+1-j}(M))_+\right).
\end{equation}
The minimum is attained. At each Fourier frequency, the positive and negative parts may each be represented by a matrix of rank at most $q$.
\end{lemma}

\begin{proof}
We first consider decompositions into positive matrices of rank at most $q$. If $P\ge0$ has eigenvalues $\mu_1\ge\cdots\ge\mu_m\ge0$, then
\begin{equation}\label{eq4.37}
 \min\left\{\sum_{k=1}^{m}k\Tr X_k:
 P=\sum_{k=1}^{m}X_k,\ X_k\ge0,\ \rank X_k\le q\right\}
 =\sum_{j=1}^m\left\lceil\frac jq\right\rceil\mu_j.
\end{equation}
Zero terms are allowed. For $P_K=\sum_{k\le K}X_k$, one has $0\le P_K\le P$ and $\rank P_K\le Kq$. Ky Fan's principle gives
\[
 \sum_{k>K}\Tr X_k=\Tr(P-P_K)\ge\sum_{j>Kq}\mu_j.
\]
Summing over $K=0,\ldots,m-1$ yields the lower bound, because the term $\mu_j$ occurs exactly $\lceil j/q\rceil$ times. Grouping the spectral vectors numbered $(k-1)q+1,\ldots,kq$ into $X_k$ attains equality.

Let $S_1+\cdots+S_q=M$. For each $a$, choose the finite optimal spectral decomposition supplied by \Cref{Lem4.4}:
\[
 S_a=\sum_{k=1}^{m}|A_{a,k}^+\rangle\langle A_{a,k}^+|
  -\sum_{k=1}^{m}|A_{a,k}^-\rangle\langle A_{a,k}^-|,
\]
where zero terms are inserted as needed. Put
\[
 X_k^\pm=\sum_{a=1}^q
 |A_{a,k}^\pm\rangle\langle A_{a,k}^\pm|,
 \qquad
 P=\sum_{k=1}^mX_k^+,
 \qquad
 Q=\sum_{k=1}^mX_k^-.
\]
Then $\rank X_k^\pm\le q$, $M=P-Q$, $P\ge M$, and $Q\ge-M$. If $p$ and $\ell$ are the positive and negative inertia indices of $M$, Weyl monotonicity gives
\[
 \lambda_j(P)\ge\lambda_j(M)>0\quad(1\le j\le p),
\]
and
\[
 \lambda_j(Q)\ge-\lambda_{m+1-j}(M)>0
 \quad(1\le j\le\ell).
\]
Applying \eqref{eq4.37} to $P$ and $Q$ proves the lower bound in \eqref{eq4.36}.

For equality, list the positive eigenvalues of $M$ as $\alpha_1\ge\cdots\ge\alpha_p>0$, with eigenvectors $u_j^+$, and the absolute values of the negative eigenvalues as $\beta_1\ge\cdots\ge\beta_\ell>0$, with eigenvectors $u_j^-$. Assign the $j$-th term of each list to the index $a=1+((j-1)\bmod q)$:
\[
 S_a=
 \sum_{\substack{\nu\ge0\\a+\nu q\le p}}
 \alpha_{a+\nu q}|u_{a+\nu q}^+\rangle
  \langle u_{a+\nu q}^+|
 -\sum_{\substack{\nu\ge0\\a+\nu q\le\ell}}
 \beta_{a+\nu q}|u_{a+\nu q}^-\rangle
  \langle u_{a+\nu q}^-|.
\]
For $j=a+\nu q$, this eigenvalue is assigned to $S_a$ at Fourier index $\nu+1=\lceil j/q\rceil$. The positive and negative terms use the corresponding positive and negative frequencies independently. The resulting cost is the expression on the right of \eqref{eq4.36}.
\end{proof}

\begin{proposition}
\label{Prop4.10}
Let $\mathbf n=(3,1,1)$, and let $T_0$ be the full block scalar torus. Suppose that $A\in\Herm(3)$ is positive semidefinite, with spectral decomposition
\[
 A=\sum_{j=1}^3\alpha_j|u_j\rangle\langle u_j|,
 \qquad \alpha_1\ge\alpha_2\ge\alpha_3\ge0,
\]
and set
\[
 H_A=\diag\bigl(A,-(\alpha_1+\alpha_3),-\alpha_2\bigr).
\]
Then $H_A\in\mathcal V_{\mathbf n,T_0}^{\mathrm H}$ and
\begin{equation*}
 \mathfrak F_{\mathbf n,T_0}(H_A)
 =d_0(e,(0,iH_A))^2
 =4\pi(\alpha_1+\alpha_2+2\alpha_3).
\end{equation*}
\end{proposition}

\begin{proof}
Since $T_0$ is the full block scalar torus, $\pi_{\mathbf n,T_0}$ is the identity on $\mathcal W_{\mathbf n}$. Write the three edge matrices as $S_{12},S_{13}\in\Herm(3)$ and $S_{23}=x\in\R$. Their images in the three blocks are
\begin{equation}\label{eq4.38}
 \mathscr A_1=S_{12}+S_{13},\qquad
 \mathscr A_2=-\Tr S_{12}+x,\qquad
 \mathscr A_3=-\Tr S_{13}-x.
\end{equation}
The first constraint therefore gives $S_{12}+S_{13}=A$. Hence \Cref{Lem4.9} and $\Phi_1(x)=4\pi|x|$ imply
\[
 \Phi_3(S_{12})+\Phi_3(S_{13})+\Phi_1(x)
 \ge\Psi_{3,2}(A)
 =4\pi(\alpha_1+\alpha_2+2\alpha_3).
\]
Equality is attained by
\[
 S_{12}
 =\alpha_1|u_1\rangle\langle u_1|
  +\alpha_3|u_3\rangle\langle u_3|,
 \qquad
 S_{13}=\alpha_2|u_2\rangle\langle u_2|,
 \qquad S_{23}=0.
\]
By \eqref{eq4.38}, its three block components are
\[
 A,\qquad-(\alpha_1+\alpha_3),\qquad-\alpha_2,
\]
so its image under $\mathscr A_{\mathbf n}$ is $H_A$. Moreover, $\Phi_3(S_{12})=4\pi(\alpha_1+2\alpha_3)$ and $\Phi_3(S_{13})=4\pi\alpha_2$. The distance identity follows from \Cref{Thm4.6}.
\end{proof}

\begin{proposition}
\label{Prop4.11}
Let $\mathbf n=(m,\underbrace{1,\ldots,1}_{q\ \mathrm{times}})$, let $T_0=\{e\}$, and let
\[
 \widehat H_A=\diag(A,0,\ldots,0),
 \qquad A\in\Herm(m),\qquad\Tr A=0.
\]
If $\lambda_1(A)\ge\cdots\ge\lambda_m(A)$ and $r_m=\lfloor m/2\rfloor$, then
\begin{align}
 \mathfrak F_{\mathbf n,\{e\}}(\widehat H_A)
 &=4\pi\min_{t\in\R}\sum_{j=1}^m
 \left\lceil\frac jq\right\rceil
 \left((\lambda_j(A)+t)_+
  +(-\lambda_{m+1-j}(A)-t)_+\right).
 \label{eq4.39}\\
 &=4\pi\sum_{j=1}^{r_m}
 \left\lceil\frac jq\right\rceil
 \bigl(\lambda_j(A)-\lambda_{m+1-j}(A)\bigr).
 \label{eq4.40}
\end{align}
All edge matrices between two scalar auxiliary blocks vanish in every minimizer. Moreover, the displayed coefficient equals
\[
 d_0(e,(0,i\widehat H_A))^2=\mathfrak F_{G,L}(i\widehat H_A).
\]

For every $q$, the minimizing shifts in \eqref{eq4.39} are
\[
 \begin{cases}
 \{-\lambda_{r_m+1}(A)\},&m=2r_m+1,\\
 [-\lambda_{r_m}(A),-\lambda_{r_m+1}(A)],&m=2r_m.
 \end{cases}
\]
If $q\ge r_m$, then all weights in \eqref{eq4.40} equal one and
\begin{equation}\label{eq4.41}
 \mathfrak F_{\mathbf n,\{e\}}(\widehat H_A)
 =4\pi\min_{t\in\R}\Tr|A+tI_m|
 =4\pi\sum_{j=1}^{r_m}
 \bigl(\lambda_j(A)-\lambda_{m+1-j}(A)\bigr).
\end{equation}
For $q=1$, \eqref{eq4.39} reduces to \eqref{eq4.33}.
\end{proposition}

\begin{proof}
For the trivial torus, orthogonal projection onto the Hermitian vertical space retains only the traceless part of the $m$-dimensional block:
\[
 \pi_{\mathbf n,\{e\}}\diag(U,u_1,\ldots,u_q)
 =\diag\left(U-\frac{\Tr U}{m}I_m,0,\ldots,0\right).
\]
Let $S_a\in\Herm(m)$ be the edge matrix joining the matrix block to the $a$-th scalar block, and let $x_{ab}\in\R$ be the edge variable joining two scalar blocks. Then
\[
 \mathscr Q_{\mathbf n,\{e\}}(S,x)
 =\diag\left(
 \sum_{a=1}^qS_a-
 \frac1m\Tr\left(\sum_{a=1}^qS_a\right)I_m,
 0,\ldots,0\right).
\]
The variables $x_{ab}$ do not enter the constraint and have nonnegative cost $4\pi\sum_{a<b}|x_{ab}|$. Hence every minimizer has $x_{ab}=0$. Feasibility is equivalent to $\sum_{a=1}^qS_a=A+tI_m$ for some $t\in\R$. Applying \Cref{Lem4.9} proves \eqref{eq4.39}.

Set $w_j=\lceil j/q\rceil$ and $W_k=\sum_{j=1}^kw_j$, with $W_0=0$. After reindexing the negative part, the function minimized in \eqref{eq4.39} is
\[
 f(t)=\sum_{j=1}^m
 \left[w_j(\lambda_j+t)_+
  +w_{m+1-j}(-\lambda_j-t)_+\right].
\]
On each open interval between consecutive breakpoints $-\lambda_j(A)$, if precisely the first $k$ shifted eigenvalues are positive, then
\[
 f'(t)=W_k-W_{m-k}.
\]
Since $W_k$ is strictly increasing, this derivative is negative for $k<m/2$, positive for $k>m/2$, and zero for $k=m/2$. At a breakpoint the subdifferential is the interval between the adjacent slopes. Hence the minimizing set is the stated median set. Substitution of a minimizing shift and pairing the upper and lower eigenvalues gives \eqref{eq4.40}.

If $q\ge r_m$, the weights occurring in the paired formula are all one, which proves \eqref{eq4.41}. If $q=1$, the weights are $j$, exactly as in the definition of $\Phi_m$.
\end{proof}

\begin{remark}
\label{Rem4.12}
After identifying $\mathfrak{pu}(m)$ with $\mathfrak{su}(m)$, the scalar term $tI_m$ represents the freedom in choosing a lift of the projective phase. Formula \eqref{eq4.40} reduces to the trace norm expression when $q\ge\lfloor m/2\rfloor$. The stronger condition $q\ge m$ occurs in the global constructions on Stiefel and Grassmann manifolds of Tanimura--Nakahara--Hayashi \cite{TanimuraNakaharaHayashi2005} and in the corresponding results of S\"onnerborn \cite{Sonnerborn2024,SonnerbornProjective2025,Sonnerborn2025}. For $\overline\Gamma_\varepsilon=[e^{i\varepsilon A}]$, S\"onnerborn's projective bound satisfies
\[
 \mathcal L_{\mathrm{proj}}(\overline\Gamma_\varepsilon)^2
 =2\pi\varepsilon\min_{t\in\R}\Tr|A+tI_m|+O(\varepsilon^2)
\]
as $\varepsilon\downarrow0$. His Grassmann length assigns squared speed $\frac12\Tr(\dot P^2)=\|Z\|_{\mathrm{HS}}^2$ to an off-diagonal velocity. The variational action of Tanimura--Nakahara--Hayashi uses the same normalization, whereas our horizontal metric gives $-\Tr(V^2)=2\|Z\|_{\mathrm{HS}}^2$. Hence \eqref{eq4.41} is twice the leading projective coefficient in this range. Formula \eqref{eq4.40} already has the trace norm form for $\lfloor m/2\rfloor\le q<m$. This is a local statement. The corresponding global criterion also depends on the spectrum of the chosen lift. If the eigenvalue $1$ has multiplicity $k_0$, the criterion in \cite{Sonnerborn2024} applies for $q\ge m-k_0$. Thus $q\ge m$ is a sufficient condition independent of the spectrum.
\end{remark}

\subsection{Distance and heat kernel consequences}
\label{Subsec4.4}

The comparison theorems quantify the nilpotent approximation in the privileged exponential chart. On every compact $\mathcal K\subset\pLie\times\lLie$, \Cref{Thm3.2} implies the uniform distance estimate
\[
 \left|\frac1r d_{\sr}\!\left(e,\exp(rX+r^2Z)\right)
 -d_0(e,(X,Z))\right|\le C_{\mathcal K}r^{1/2}.
\]
Under the symmetric hypotheses of \Cref{Thm3.4}, the same statement holds on compact subsets of $\pLie\times\kLie$, with remainder $O_{\mathcal K}(r)$. In the compact case this distance estimate also has order $O_{\mathcal K}(r)$ on sets separated from $(0,0)$. The squared estimate remains uniform even when the leading distance vanishes.

The one-column AIII coefficient was evaluated explicitly in \cite[Theorem~6.2]{BauerHanYang2026}. We now apply the comparison theorems to the general block coefficient. In the two-block full Levi case, the symmetric comparison improves the error in absolute value from $O_R(\varepsilon^{3/2})$ to $O_R(\varepsilon^2)$.

\begin{corollary}
\label{Cor4.13}
For every $R>0$, uniformly for $\mathfrak F_{\mathbf n,T_0}(H)\le R$,
\begin{equation*}
 d_{\sr}(e,\exp(i\varepsilon H))^2
 =\varepsilon\mathfrak F_{\mathbf n,T_0}(H)
  +O_R(\varepsilon^{3/2})
\end{equation*}
as $\varepsilon\downarrow0$. If $s=2$ and $T_0$ is the full block scalar torus, then the pair is symmetric and the remainder improves to $O_R(\varepsilon^2)$.
\end{corollary}

\begin{proof}
The set $\{H:\mathfrak F_{\mathbf n,T_0}(H)\le R\}$ is compact by \eqref{eq4.24}. Apply \Cref{Thm3.2} with $X=0$, $Z=iH$, and $r=\sqrt\varepsilon$, then use \eqref{eq4.23}. For two blocks and $T_0=Z_{\mathbf n}$, the horizontal space is precisely the off-diagonal part of the AIII symmetric decomposition, so \Cref{Thm3.4} applies.
\end{proof}

The same expansions determine the logarithmic heat kernel rate in an iterated limit. The heat time first tends to zero, and the endpoint then tends to the identity. The groups considered here are unimodular. This is automatic in the compact case. Under the symmetric hypotheses, $\gLie=\pLie+[\pLie,\pLie]$ and $[\pLie,\pLie]\subset\kLie$ imply $\kLie=[\pLie,\pLie]$. For $X\in\pLie$, the map $\ad_X$ is off diagonal with respect to $\gLie=\kLie\oplus\pLie$, and hence has trace zero. Every element of $\kLie$ is a linear combination of brackets $[X,Y]$, and
\[
 \Tr(\ad_{[X,Y]})=\Tr[\ad_X,\ad_Y]=0.
\]
Thus the modular character vanishes, and the left-invariant horizontal fields have zero divergence with respect to Haar measure.

\begin{corollary}
\label{Cor4.14}
Let $X_1,\ldots,X_d$ be an orthonormal basis of $\pLie$ and let $\widetilde X_j$ denote the corresponding left-invariant vector fields. Set
\[
 \Delta_{\pLie}=\sum_{j=1}^d\widetilde X_j^{\,2}.
\]
Let $\mu$ be a Haar measure on $G$, normalized to have total mass one when $G$ is compact. Write $p_t$ for the convolution heat kernel, so that the integral kernel of $e^{t\Delta_{\pLie}}$ with respect to $\mu$ is
\[
 K_t(g,h)=p_t(g^{-1}h).
\]
For fixed $r>0$, define
\begin{equation}\label{eq4.42}
 \mathcal I_r(X,Z)
 := -\lim_{t\downarrow0}4t
 \log p_t\!\left(\exp(rX+r^2Z)\right).
\end{equation}
The limit in \eqref{eq4.42} is taken before $r\downarrow0$. Under the compact hypotheses of \Cref{Thm3.2}, every compact $\mathcal K\subset\pLie\times\lLie$ satisfies
\begin{equation*}
 \mathcal I_r(X,Z)=r^2d_0(e,(X,Z))^2+O_{\mathcal K}(r^3)
\end{equation*}
uniformly for $(X,Z)\in\mathcal K$. Under the symmetric hypotheses of \Cref{Thm3.4}, the same conclusion holds for compact $\mathcal K\subset\pLie\times\kLie$, with remainder $O_{\mathcal K}(r^4)$. At a datum satisfying the conditions in \Cref{Ass3.6},
\[
 \mathcal I_r(X,Z)
 =r^2d_0(e,(X,Z))^2+r^4 c_4(X,Z)+O(r^6),
 \qquad
 c_4(X,Z)=-\lambda_0\bigl(Q(u_0),R(u_0)\bigr).
\]
At a block vertical endpoint, $d_0(e,(0,iH))^2=\mathfrak F_{\mathbf n,T_0}(H)$.
\end{corollary}

\begin{proof}
For the present generator, L\'eandre's logarithmic theorem \cite{Leandre1987Upper,Leandre1987Lower} gives
\[
 -\lim_{t\downarrow0}4t\log p_t(g)
 =d_{\sr}(e,g)^2.
\]
Thus $\mathcal I_r(X,Z)=d_{\sr}(e,\exp(rX+r^2Z))^2$ for every fixed $r$. The stated expansions now follow from \Cref{Thm3.2,Thm3.4} and \Cref{Thm3.7,Thm4.6}.
\end{proof}

\section*{Acknowledgments}

\noindent\textbf{Funding.} This work was supported by the National Natural Science Foundation of China (12301145, 12561020, 12261107) and the Yunnan Fundamental Research Projects (202401AU070123, 202601AT070048).

\medskip
{\bf Author Contributions:} All authors contributed equally to the writing and preparation of the manuscript.

\medskip
{\bf Data availability:} Data sharing is not applicable to this article as no new data were created or analyzed in this study.

\medskip
{\bf Conflict of Interests:} The authors declare that they have no conflict of interest.

\medskip
{\bf AI Declaration:} The authors used GPT-6 (OpenAI) to assist with literature searches and English language editing. All mathematical arguments and proofs in the final manuscript were written and checked by the authors, who take full responsibility for its content.

\end{document}